
\input amstex
\documentstyle{amsppt}

\magnification=\magstep1
\hsize=6.5truein
\vsize=9truein


\font \smallrm=cmr10 at 10truept
\font \smallbf=cmbx10 at 10truept
\font \smallit=cmti10 at 10truept
 at 10truept
 at 10truept

\baselineskip=.15truein

\def \strangearrow {\,\succ\joinrel\!\!\relbar\joinrel\rightsquigarrow\,}

\def \unon {1,\dots,n}

\def \E {\Bbb E \,}
\def \N {\Bbb N}
\def \Z {\Bbb Z}
\def \Q {\Bbb Q}
\def \C {\Bbb C}
 \def \Lam {\Lambda \hskip-6,1pt \Lambda}

\def \tm {t^{-1}}
\def \qm {q^{-1}}
\def \Cq {\C(q)}

\def \Reps {R_\varepsilon}

\def \gerg {{\frak g}}
\def \gerh {{\frak h}}

\def \gerU {{\frak U}}

\def \gerB {{\frak B}}
\def \ghat {{\hat \gerg}}
\def \hhat {{\hat \gerh}}

\def \calF {\hbox{$ {\Cal F} $}}

\def \edot {\dot E}
\def \fdot {\dot F}
\def \etilde {\tilde E}
\def \ftilde {\tilde F}
\def \ehat {\hat E}
\def \fhat {\hat F}
\def \echeck {\check E}
\def \fcheck {\check F}

\def \Edot {\dot \E}
\def \Etilde {\tilde \E}
\def \Ehat {\hat \E}

\def \phire {\Phi^{\text{re}}}
\def \phiim {\Phi^{\text{im}}}
\def \phipre {\Phi_+^{\text{re}}}
\def \phipim {\Phi_+^{\text{im}}}
\def \phitildep {\widetilde{\Phi}_+}
\def \phitildepim {\widetilde{\Phi}_+^{\text{im}}}

\def \uqg {U_q(\ghat)}
\def \gerUqg {\gerU_q(\ghat)}

\def \calUqg {{\Cal U}_q(\ghat)}

\def \uqm {U_q^-}
\def \gerUqm {\gerU_q^-}
\def \calUqm {{\Cal U}_q^-}
\def \uqz {U_q^0}
\def \calUqz {{\Cal U}_q^0}
\def \gerUqz {\gerU_q^0}
\def \uqp {U_q^+}
\def \gerUqp {\gerU_q^+}
\def \calUqp {{\Cal U}_q^+}

\def \uqbm {U_q^\leq}
\def \calUqbm {{\Cal U}_q^\leq}
\def \gerUqbm {\gerU_q^\leq}
\def \uqbp {U_q^\geq}
\def \calUqbp {{\Cal U}_q^\geq}
\def \gerUqbp {\gerU_q^\geq}

\document

\topmatter

{\ }

\vskip-33pt

\hfill   {{\sl Communications in Algebra\/}  {\bf 27},
no.~2 (1999), 903--918}
\hskip19pt   {\ }

\vskip41pt

 \title
   A PBW basis for Lusztig's form of untwisted affine quantum groups
 \endtitle

\author
  Fabio Gavarini
\endauthor

\affil
   Universit\`a degli Studi di Roma ``Tor Vergata'' --- Dipartimento
di Matematica  \\
   Via della Ricerca Scientifica, I-00133 Roma --- ITALY  \\
\endaffil

\address\hskip-\parindent
        Universit\`a degli Studi di Roma ``Tor Vergata''  \newline
        Dipartimento di Matematica  \newline
        Via della Ricerca Scientifica  \newline
        I-00133 Roma --- ITALY  \newline
        e-mail: \  gavarini\@axp.mat.utovrm.it,  \newline
        \phantom{e-mail:} \  gavarini\@mat.uniroma2.it,  \newline
        \phantom{e-mail:} \  gavarini\@mat.uniroma3.it
\endaddress

 \abstract
   Let  $ \ghat $  be an untwisted affine Kac-Moody algebra over the field
$ \C $,  and let  $ \uqg $  be the associated quantum enveloping algebra;
let  $ \gerUqg $  be the Lusztig's integer form of  $ \uqg $,  generated
by  $ q $-divided powers of Chevalley generators over a suitable subring
$ R $  of  $ \Cq $.  We prove a Poincar\'e-Birkhoff-Witt like theorem for
$ \gerUqg $,  yielding a basis over  $ R $  made of ordered products of
$ q $--divided  powers of suitable quantum root vectors.
 \endabstract

\endtopmatter

\footnote""{ Keywords:  {\sl PBW Theorems, Restricted Affine Quantum Groups} }


\footnote""{ 1991 {\it Mathematics Subject Classification,}
Primary 17B37, 81R50 }
\footnote""{ Partially supported by a post-doc fellowship of the
{\it Consiglio Nazionale delle Ricerche} \, (Italy) }

\vskip0,485truecm

 \centerline{ \bf  Introduction }

\vskip10pt

\hfill  \hbox{\vbox{ \hbox{\it \hskip9pt   "Questa forma \`e duale }
                     \hbox{\it \hskip19pt   di un'altra gi\`a nota }
                     \hbox{\it             che ha un suo teorem PBW. }
                     \hbox{\it \hskip22pt     Ed \`e subito base" }
                     \vskip4pt
                     \hbox{\sl \hskip11pt    N.~Barbecue, "Scholia" } }
\hskip1truecm }

\vskip8pt
   Let  $ \ghat $  be an untwisted affine Kac-Moody algebra, and let
$ \uqg $  be the associated quantum enveloping algebra.  In [Be1], [Be2],
quantum root vectors are defined, and a basis of Poincar\'e{}-Birkhoff-Witt
type for  $ \uqg $  is constructed, made of ordered monomials in the quantum
root vectors.
                                                   \par
  Now let  $ \gerUqg $  be the Lusztig's integer form of  $ \uqg $,
generated over  $ \Z \left[ q, \qm \right] $  by  $ q $--divided  powers
$ E_i^{(n)} $,  $ F_i^{(m)} $;  for technical reasons, we shall use a larger
ground ring  $ R $.  In this paper we find a PBW basis of  $ \gerUqg $
as an $ R $--module,  made of ordered products of  $ q $--divided  powers
of (suitable renormalizations of) quantum root vectors.
                                                   \par
  As a first step we reduce the problem to finding a basis for  $ \gerUqp $
(the positive part of  $ \gerUqg $).  Second, we exploit the duality among
PBW basis in  $ \uqp $  and  in  $ \uqm $  --- proved in [Da2] ---  to get
from there our key result, namely finding a basis of  $ \gerUqp $.
                                                   \par
  Such an approach is entirely different from the classical ones, to be found
in [Ga] and [Mi]; on the other hand, the comparison with the classical setting
is quite interesting: this is sketched in the last section, where also a second
PBW theorem is proved and some further conjectures are presented.

\vskip7pt

\centerline { ACKNOWLEDGEMENTS }

\vskip4pt

  The author wishes to thank J.~Beck and V.~Chari for explanations about their
works, as well as I.~Damiani and M.~Rosso for many helpful discussions.

\vskip1,7truecm

   \centerline{ \bf  \S \; 1 \,  Notations }

\vskip10pt

   {\bf 1.1  The classical data.}  \  We shall adopt notations of [Da2],
which we recall in this section (they follow [Bo] and [Ka]).
                                                          \par
  Let  $ \gerg $  be a simple finite dimensional Lie algebra over the field
$ \C $  of complex numbers, and consider the folllowing data.
                                                          \par
   $ I_0 = \{1, \dots, n \} \, $,  the set of vertices of the Dynkin
diagram of  $ \gerg $  (see [Da1] for the identification between  $ I_0 $
and  $ \{1, \dots, n \} $);  $ \gerh $  a Cartan subalgebra of  $ \gerg $,
generated by   $ \, \{h_1, \dots, h_n\} = \{\, h_i \mid i \in I_0 \,\} \, $;
$ \Phi_0 = \Phi_{0,+} \cup \left( -\Phi_{0,+} \right) \subseteq \gerh^* \, $
the root system of  $ \gerg $,  with  $ \, \Phi_{0,+} \, $  the set of positive
roots, and  $ \, \Pi := \{\alpha_1, \dots, \alpha_n\} = \{\, \alpha_i \mid i
\in I_0 \,\} \, $  the set of simple roots;  $ Q_0 := \sum_{\alpha \in \Phi_0}
\Z \alpha = \oplus_{i \in I_0} \Z \alpha_i \, $  the root lattice of
$ \gerg $,  and  $ Q_0^\vee := \sum_{\alpha \in \Phi_0} \Z \alpha^\vee =
\oplus_{i \in I_0} \Z \alpha_i^\vee \, $  the coroot lattice;  $ W_0 $  the
Weyl group of  $ \gerg $.  Finally, we fix a function  $ \, o : I_0
\longrightarrow \{ \pm 1 \} \, $  such that  $ \, a_{ij} < 0 \Longrightarrow
o(i) \, o(j) = -1 \, $.
                                                          \par
  We denote  $ \ghat $  the untwisted affine Kac-Moody algebra associated to
$ \gerg $  and we consider its loop-algebra like realization as
  $$  \ghat = \gerg \otimes_\C \C \left[ t, \tm \right] \oplus \C \cdot
c \oplus \C \cdot \partial  $$
with the Lie bracket given by:  $ \; [c,z] = 0 \, $,  $ \, \left[ \partial, x
\otimes t^m \right] = m x \otimes t^m \, $,  $ \, \left[ x \otimes t^r, y
\otimes t^s \right] = [x,y] \otimes t^{r+s} + \delta_{r,-s} r \, (x \vert y) \,
c \; $  for all  $ \, z \in \ghat $,  $ x, y \in \gerg $,  $ m \in \Z $,  where
$ \, (\,\cdot\, \vert \,\cdot\, ) \, $  is the Killing form of  $ \gerg $,
normalized in such a way that  $ \, \big( h_i, h_j \big) = {\, a_{ij} \, \over
\, d_j \,} \, $.
                                                          \par
  For  $ \ghat $  we define:  $ \, I := \{0,1,\dots,n\} \supset I_0 \, $  to be
the set of vertices of the Dynkin diagram, and  $ \, I_\infty := I \, \cup
\{\infty\} \, $;  $ \, A = {(a_{ij})}_{i,j \in I} \, $  the (generalized)
Cartan matrix and  $ \, D = \text{diag}(d_0, d_1, \dots, d_n) \, $  the
diagonal matrix with relatively prime positive integral entries such that  $ \,
D A \, $  is symmetric  ($ d_0 = 1 $);  $ \hhat := \gerh \oplus \C \cdot c
\oplus \C \cdot \partial \, (\, \subseteq \ghat \,) \, $;  $ \, \Phi = \Phi_+
\cup (-\Phi_+) \, \left(\, \subset {(\gerh \oplus \C \cdot c)}^* \subset
\hhat^* \, \right) \, $  the root system,  $ \, \Phi_+ = \phipre \cup \phipim
\, $  the set of positive roots,  $ \, \{ \alpha_0, \alpha_1, \dots, \alpha_n
\} = \{\, \alpha_i \mid i \in I \,\} \, $  the set of simple roots,  $ \,
\phipim = \, \{\, m \delta \mid m \in \N_+ \, \} \, $  the set of imaginary
positive roots (where  $ \, \delta = \sum_{i \in I} d_i \alpha_i = \theta +
\alpha_0 \, $  and  $ \theta $  is the longest positive root of  $ \gerg $),
$ \, \phipre = \Phi_{0,+} \cup \{\, \alpha + m \delta \mid \alpha \in \Phi_0 \,
,  m > 0 \,\} \, $  the set of real positive roots.  Then  $ \ghat $  has a
decomposition into direct sum of  $ \hhat $  and root spaces  $ \, \ghat =
\hhat \oplus \big( \oplus_{\alpha \in \Phi} \ghat_\alpha \big) \, $,  and
  $$  \text{dim}_{\Bbb C} \big( \ghat_\alpha \big) = 1  \;\;\; \forall\; \alpha
\in \phire \, ,  \qquad  \text{dim}_{\Bbb C} \big( \gerh \big) = \#
\big(I_0\big) = n \;\;\; \forall\; \alpha \in \phiim \; ;  $$
therefore we define the set  $ \phitildep $  of "positive roots with
multiplicity" as  $ \, \phitildep := \phipre \cup \phitildepim \, $,  where
$ \, \phitildepim := \phipim \times I_0 \, $;  then we denote  $ \, p :
\phitildep \rightarrow \Phi_+ \, $  the natural projection map.  Furthermore,
we have:  the root lattice (of  $ \ghat $)  $ \, Q = \sum_{\alpha \in \Phi} \Z
\cdot \alpha = \oplus_{i \in I} \Z \cdot \alpha_i = \Z \cdot \alpha_0 \oplus Q_0 = Q_0 \oplus \Z \cdot \delta \, $,  $ \, Q_\infty := Q \oplus \Z \cdot
\alpha_\infty \, $  (where  $ \, \alpha_\infty (\in \gerh^*) \, $  is the dual
element to  $ \, \delta (\in \gerh)) $,  and the order relation  $ \, \leq
\, $  on  $ Q_\infty $  given by  $ \, \alpha \leq \beta \iff \beta - \alpha \in
Q_+ \, $,  with  $ \, Q_+ := \sum_{i \in I} \N \cdot \alpha_i \, $;  the
non-degenerate symmetric bilinear form on  $ \, Q_\infty \otimes_\Z {\Bbb R}
\, $  given by  $ \, (\alpha_i \vert \alpha_j) = d_i a_{ij} \, (\, \forall\, i,
j \in I \,) $,  $ (\alpha_\infty \vert Q_0 \oplus \Z \cdot \alpha_\infty) = 0 $,
$ (\alpha_\infty \vert \delta) = 1 \, $;  the group  $ \, W = W_0 \ltimes
Q_0^\vee \, $,  the subset of simple reflections  $ \, \{ s_0, s_1, \dots, s_n
\} = \{\, s_i \mid i \in I \,\} \, (\, \subseteq W \,) \, $,  and the length
function  $ \, l \colon \, W \rightarrow \N \, $;  the braid group
$ {\Cal B} $  (associated to  $ W \, $),  generated by  $ \, \{ T_0, T_1,
\dots, T_n \} = \{\, T_i \mid i \in I \,\} \, $,  and the section  $ \, T
\colon \, W \rightarrow {\Cal B} \, $  such that  $ \, T_w = T_{i_1} \cdots
T_{i_r} \, $  for all  $ \, w = s_{i_1} \cdots s_{i_r} \in W \, $  with  $ l(w)
= r $.

\vskip7pt

   {\bf 1.2  The quantum group  $ \uqg $.}  \  The quantized  universal
enveloping algebra  $ \uqg $  (cf.~e.g.~[Dr], [Lu2]) is the unital associative
$ \Cq $--algebra  with generators
   $$  F_1, \dots , F_n,  K_\nu \, (\, \nu \in Q_\infty \,),  \, E_1, \dots,
E_n  $$
and relations
  $$  \displaylines {
   \hfill   K_\mu K_\nu = K_{\mu + \nu} = K_\nu K_\mu  \quad , \qquad K_0 = 1
\hfill  \forall \, \mu, \nu \in Q_\infty  \cr
   \hfill   K_{\mu} E_i = q^{(\mu \vert \alpha_i)} E_i K_{\mu} \, ,
\quad  K_{\mu} F_i = q^{-(\mu \vert \alpha_i)} F_i K_{\mu}  \hfill  \forall \,
\mu \in Q_\infty, \; \forall \, i \in I  \cr
   \hfill   E_i F_j - F_j E_i = \delta_{ij} {K_{\alpha_i} - K_{-\alpha_i} \over
q_i - q_i^{-1}}   \hfill  \forall \, i, j \in I  \cr
   \sum_{k = 0}^{1-a_{ij}} (-1)^k {\left[ { 1-a_{ij} \atop k } \right]}_{q_i}
\!\! E_i^{1-a_{ij}-k} \! E_j E_i^k = 0 \, ,  \,  \sum_{k = 0}^{1-a_{ij}} (-1)^k
{\left[{1-a_{ij} \atop k} \right]}_{q_i} \!\! F_i^{1-a_{ij}-k} \! F_j F_i^k = 0
 \hfill  \forall \, i \ne j  \cr }  $$
where  $ \; {\left[{m \atop n} \right]}_q := {{[m]}_q ! \over {[m-n]}_q ! \,
{[n]}_q !} \, $,  $ \, {[k]}_q ! := \prod_{s=1}^k {[s]}_q \, $,  $ \, {[s]}_q :=
{q^s - q^{-s} \over q - q^{-1}} \, $,  for all  $ m $,  $ n $,  $ k $,  $ s \in
\N_+ $,  $ n \leq m $,  with
$ {[s]}_q , {[k]}_q ! , {\left[{m \atop n} \right]}_q \in \Z \left[ q, \qm
\right] $.  For later use, we define also
  $$  q_\alpha := q^{\,(\alpha \vert \alpha)\, \over \,2\,}  \; \forall \,
\alpha \in \phipre \, ,  \quad  q_\alpha := q^{d_i} \; \forall \, \alpha = (r
\delta, i) \in \phitildepim \, ,  \quad  q_i := q_{\alpha_i} = q^{d_i}  \;
\forall \, i \in I \, .  $$
   \indent   A Hopf algebra structure on  $ \, \uqg \, $  is defined by
($ i= \unon; \, \mu \in Q_\infty $)
  $$  \matrix
   \Delta(F_i) := F_i \otimes K_{-\alpha_i} + 1 \otimes F_i \,,  &
S(F_i):= -F_i K_{\alpha_i} \, ,  &  \epsilon(F_i) := 0 \, \phantom{.}  \\
   \Delta(K_\mu) := K_\mu \otimes K_\mu \, ,  &   S(K_\mu):= K_{-\mu} \, ,
&  \epsilon(K_\mu) := 1 \, \phantom{.}  \\
  \Delta(E_i) := E_i \otimes 1 + K_{\alpha_i} \otimes E_i \, ,  &
S(E_i):= -K_{-\alpha_i} E_i \, ,  &  \epsilon(E_i) := 0 \, .  \\
      \endmatrix  $$
   \indent   Moreover,  $ \uqg $  has a natural algebra grading  $ \, \uqg =
\oplus_{\eta \in Q} {\uqg}_\eta \, $.
                                                \par
  We have a  $ \C $--antilinear  antiinvolution  $ \, \Omega : \uqg \rightarrow
\uqg \, $  defined by
  $$  \Omega(q) := \qm \, ,  \quad  \Omega(E_i) := F_i \, ,  \quad
\Omega(K_\nu) := K_{-\nu} \, ,  \quad  \Omega(E_i) := F_i \, ,  \qquad
\forall \, i \in I, \nu \in Q_\infty  $$
and a braid group action on  $ \uqg $  wich commutes with  $ \Omega $.
                                                \par
  Let  $ \, \uqp $,  $ \uqz $,  $ \uqm $  be the subalgebras of  $ \uqg $
respectively generated by  $ \, \{\, E_i \mid i \in I \,\} $,  $ \{\, K_\nu
\mid \nu \in Q_\infty \,\} $,  $ \{\, F_i \mid i \in I \,\} $;  let also
$ \, \uqbp := \uqp \cdot \uqz = \uqz \cdot \uqp \, $,  $ \, \uqbm := \uqm
\cdot \uqz = \uqz \cdot \uqm \, $,  to be called  {\it quantum Borel
(sub)algebras}:  these are  {\sl Hopf}  subalgebras of  $ \uqg $.  Finally,
multiplication provides linear isomorphisms ("triangular decompositions")
  $$  \displaylines{
   \uqg \cong \uqp \otimes \uqz \otimes \uqm \cong \uqm \otimes \uqz \otimes
\uqp  \cr
   \uqbp \cong \uqp \otimes \uqz \cong \uqz \otimes \uqp \, ,  \quad  \uqbm
\cong \uqz \otimes \uqm \cong \uqm \otimes \uqz \, .  \cr }  $$

\vskip4pt

   {\it Remark:} \, In the definition of  $ \uqg $  several choices for the
"toral part"  $ \uqz $  are possible, mainly depending on the choice of any
lattice  $ M $  such that  $ \, Q \leq M \leq P \, $,  $ P $  being the
weight lattice of  $ \gerg $  (cf.~for instance [B-K]).  All what follows
holds as well for every such choice, up to suitably adapting the statements
involving the toral part.

\vskip7pt

   {\bf 1.3  The DRT pairing.}  \  Let  $ {\left( {\uqbm}
\right)}^{\text{op}} $  denote the Hopf algebra with the same structure of
$ \uqbm $  but for the coproduct, which  is turned into the opposite one.
Then a perfect (= non-degenerate) pairing of Hopf algebras  $ \, \pi : \uqbp
\times {\left({\uqbm}\right)}^{\text{op}} \longrightarrow \Cq \, $  exists,
defined by
  $$  \pi (K_\lambda, K_\mu) := q^{-(\lambda \vert \mu)} \, ,  \quad  \pi
(K_\lambda, F_i) := 0 =: \pi (E_i, K_\mu) \, ,  \quad  \pi (E_i,F_j) := {\,
\delta_{ij} \, \over \, q_i^{-1} - q_i \,}  $$
for all  $ \, i, j \in I $,  $ \lambda, \mu \in Q_\infty \, $  (cf.~[Ta]).  This
is a graded pairing, i.~e.  $ \, \pi{\big\vert}_{{\left( \uqbp \right)}_\eta
\times {\left( {\left( {\uqbm} \right)}^{\text{op}} \right)}_\zeta} = 0 \, $
for all  $ \, \eta, \zeta \in Q_+ \, $  such that  $ \, \eta + \zeta \neq 0
\, $  (the grading is the one inherited from  $ \uqg $),  and so  $ \,
\pi{\big\vert}_{{\left( \uqbp \right)}_\eta \times {\left( {\left( {\uqbm}
\right)}^{\text{op}} \right)}_{-\eta}} \, $  is non-degenerate, for all  $ \,
\eta \in Q_+ \, $;  finally, we have
  $$  \pi (x K_\lambda, y) = \pi (x,y) = \pi(x, y K_\mu) \quad \forall \, x
\in \uqbp \, ,  \quad \forall \, y \in \uqbm \, ,  \; \forall \, \lambda, \mu
\in Q_\infty \, .   \eqno (1.1)  $$

\vskip1,7truecm

   \centerline{ \bf  \S \; 2 \,  Quantum root vectors and PBW bases }

\vskip10pt

   {\bf 2.1  Ordering positive roots.}  \  In this section we sketch the
construction of PBW bases, provided by Beck (cf.~[Be1], [Be2]), following [Da2],
\S 2.  We can take a suitable function  $ \, \iota : \Z \rightarrow I \, $
(fixed once and for all) such that, letting
   $$  \beta_k \equiv \beta_k^{(\iota)}
    := \cases
          s_{\iota(1)} s_{\iota(2)} \cdots s_{\iota(k-1)} \left(
\alpha_{\iota(k)} \right)  &  \quad  \text{for all}  \  k \geq 1  \\
          s_{\iota(0)} s_{\iota(-1)} \cdots s_{\iota(k+1)} \left(
\alpha_{\iota(k)} \right)  &  \quad  \text{for all}  \  k \leq 0  \\
       \endcases  $$
then  $ \, k \mapsto \beta_k \, $  defines a bijection of  $ \Z $  onto
$ \phipre \, $,  so that  $ \, \big\{\, \beta_k \,\big\vert\, k \geq 1 \,\big\}
= \big\{\, r \delta - \alpha \,\big\vert\, r > 0, \alpha \in \Phi_{0,+}
\,\big\} \, $  and  $ \, \big\{\, \beta_k \,\big\vert\, k \leq 0 \,\big\} =
\big\{\, r \delta + \alpha \,\big\vert\, r \geq 0, \alpha \in \Phi_{0,+}
\,\big\} \, $.
                                          \par
  Then one defines a total ordering on  $ \phitildep $ as follows:
  $$  \displaylines{
   {\ }  \beta_1 \preceq \beta_2 \preceq \beta_3 \preceq \cdots \preceq
\beta_{k-1} \preceq \beta_k \preceq \beta_{k+1} \preceq \cdots  \hfill {\ }
\cr
   \cdots \preceq \big( (r+1) \delta,n) \preceq (r\delta,1) \preceq (r\delta,2)
\cdots \preceq (2\delta,n) \preceq (\delta,1) \preceq (\delta,2) \preceq \cdots
\preceq (\delta,n) \preceq \cdots  \cr
   {\ } \hfill  \cdots \preceq \beta_{-(k+1)} \preceq \beta_{-k} \preceq
\beta_{-(k-1)} \cdots \preceq \beta_{-2} \preceq \beta_{-1} \preceq \beta_0
{\ }  \cr }  $$

\vskip7pt

   {\bf 2.2  Quantum root vectors.}  \  For  $ \iota $  fixed as in \S 2.1, we
define quantum root by
  $$  E_{\beta_k}^{(\iota)} \equiv E_{\beta_k} :=
     \cases
        T_{\iota(1)} T_{\iota(2)} \cdots T_{\iota(k-1)} \left( E_{\iota(k)}
\right)  &  \quad  \text{for all}  \  k \geq 1  \\
          T_{\iota(0)}^{-1} T_{\iota(-1)}^{-1} \cdots T_{\iota(k+1)}^{-1}
\left( E_{\iota(k)} \right)  &  \quad  \text{for all}  \  k \leq 0  \\
     \endcases  $$
  $$  \displaylines{
   \etilde_{(r \delta, i)}^{(\iota)} \equiv \etilde_{(r \delta, i)} :=
q_i^{-2} E_i E_{r \delta - \alpha_i} - E_{r \delta - \alpha_i} E_i  \qquad
\quad \forall \, r > 0,  \; \forall \, i \in I_0  \cr
   \hfill   \left( q_i - q_i^{-1} \right) \sum_{r>0} E_{(r \delta,i)} \cdot
\zeta^r := \log \left( 1 - \left( q_i - q_i^{-1} \right) \sum_{s>0} \etilde_{(s
\delta,i)} \cdot \zeta^r \right)   \hfill   (2.1)  \cr
   \ftilde_{(r \delta,i)} := \Omega \left( \etilde_{(r \delta,i)} \right)
\; \quad \forall \, (r \delta,i) \in \phitildepim \, ,  \qquad  F_\alpha :=
\Omega \left( E_\alpha \right) \; \quad \forall \, \alpha \in \phitildep
\cr } $$
To be complete, we define also  $ \, \etilde_{(0,i)} = E_{(0,i)} := 1 =:
\ftilde_{(0,i)} = F_{(0,i)} \, $  for all  $ \, i \in I_0 \, $.
                                                            \par
   For later use, we define other imaginary root vectors by the
recursive formulas
  $$  \edot_{(0 \delta,i)} := 1 \, ,  \qquad  \edot_{(r \delta,i)} := q_i^r {\,
1 \, \over \, {[r]}_{q_i} \,} \sum_{s=1}^r \etilde_{(s \delta, i)}
\edot_{((r-s) \delta,i)}  \qquad \forall \, (r \delta,i) \in \phitildepim
\eqno (2.2)  $$
and  $ \, \fdot_{(r \delta,i)} := \Omega \left( \edot_{(r \delta,i)} \right)
\, $;  the relationship among the  $ \edot_{(h \delta, i)} $'s  and  the
$ E_{(k \delta, i)} $'s                          is given by\break
  $$  \edot_{(0,i)} = 1 = E_{(0,i)} \, ,  \quad  \edot_{(r \delta, i)} =
- {\, 1 \, \over \, r \,} \sum_{s=1}^r q_i^s {\, s \, \over \, {[s]}_{q_i} \,}
E_{(s \delta, i)} \edot_{((r-s) \delta, i)}  \quad \forall \, (r \delta,i) \in
\phitildepim   \eqno (2.3)  $$
(cf.~[C-P], \S 3) and similarly with  "$ \! F \, $"  instead of  "$ \! E \, $".
 Definitions give
  $$  E_{\beta} \in {\left( \uqp \right)}_\beta  \quad \forall \, \beta \in
\phipre \, ,  \qquad  \etilde_{(r \delta,i)} \, , \, E_{(r \delta,i)} \, ,
\, \edot_{(r \delta,i)} \in {\left( \uqp \right)}_{r \delta}  \quad  \forall
\,\, r \in \N \, , \, k \in N_+ \, , \, i \in I_0  $$
and similarly for negative root vectors, with  "$ \! F \, $"'s  instead of
"$ \! E \, $"'s.  Finally we recall the following property of imaginary root
vectors, which will be useful later:

\vskip3pt

   \centerline {\it  All quantum root vectors attached to imaginary roots
commute with each other.}

\vskip7pt

   {\bf 2.3  PBW bases and orthogonality.}  \  It is proved in [Be2] that the
set of ordered monomials in the root vectors  $ E_\alpha $'s  (according
to the order  $ \, \preceq \, $  on  $ \Phi_+ $  defined in \S 2.1),
namely the  $ \, \prod_{\alpha \in \Phi_+} E_\alpha^{n_\alpha} $'s
(where the  $ \, n_\alpha \in \N \, $  are almost all zero), is a
$ \Cq $--basis  of  $ \uqp $;  similarly, the set of ordered
monomials in the root vectors  $ F_\alpha $'s is a  $ \Cq $--basis  of
$ \uqm $.  Since clearly  $ \, \{\, K_\alpha \mid \alpha \in Q_\infty \,\} =
\left\{\, \prod_{i \in I_\infty} K_i^{l_i} \,\Big\vert\, l_i \in \Z \;
\forall \, i \in I_\infty \,\right\} \, $  is a  $ \Cq $--basis  of  $ \uqz $,
from triangular decompositions one concludes that the sets of ordered
monomials  $ \, \prod_{\alpha \in \Phi_+} E_\alpha^{n_\alpha} \prod_{i \in
I_\infty} K_i^{l_i} \prod_{\alpha \in \Phi_+} F_\alpha^{m_\alpha} \, $  or
$ \, \prod_{\alpha \in \Phi_+} F_\alpha^{m_\alpha} \prod_{i \in I_\infty}
K_i^{l_i} \prod_{\alpha \in \Phi_+} E_\alpha^{n_\alpha} \, $  are
$ \Cq $--bases  of  $ \uqg $,  and similarly for quantum Borel subalgebras.
From [Da2] we recall also
  $$  \displaylines{
   \pi \left( E_\alpha , F_\beta \right) = {\, \delta_{\alpha,\beta} \, \over
\, \left( q_\alpha^{-1} - q_\alpha \right) \,} \, ,  \;\;  \pi \left( E_\alpha
, F_\gamma \right) = 0 \, ,  \;\;  \pi \left( E_\gamma , F_\alpha \right) = 0
\qquad  \forall \, \alpha, \beta \in \phipre \, ,
\, \gamma \in \phitildepim  \cr
   \pi \left( E_{(r \delta, i)} , F_{(s \delta, j)} \right) = \delta_{r,s}
{\big( o(i) o(j) \big)}^r {\, {[r a_{ij}]}_{q_i} \, \over \, r \left( q_j^{-1}
- q_j \right) \,}  \qquad \forall \, (r, \delta, i), (s \delta, j) \in
\phitildepim \, ;  \cr }  $$
these formulas are the starting point to build up orthogonal bases of
$ \uqp $  and  $ \uqm $.

\vskip7pt

\proclaim{Lemma 2.4}  For all  $ r \in \N_+ \, $,  let  $ V_r $,
resp.~$ W_r $,  be the  $ \Cq $--vector  space with basis  $ \{\, E_{(r \delta,
i)} \,\vert\, i \in I_0 \,\} $, resp.~$ \{\, F_{(r \delta, i)} \,\vert\, i \in
I_0 \,\} $,  and let  $ \, \{\, x_{r,i} \mid i \in I_0 \,\} \, $  and  $ \, \{\,
y_{r,j} \mid j \in I_0 \,\} \, $  be bases of  $ V_r $  and  $ W_r $  orthogonal
of each other with respect to  $ \pi $,  namely  $ \, \pi \big( x_{r,i}, y_{r,j}
\big) = 0 \, $  for all  $ i \neq 0 $  ($ i, j \in I_0 $).  Then  $ \, \Big\{\,
\prod_{k \leq 0} E_{\beta_k}^{n_k} \cdot \! \prod_{r \in \N, i \in I_0}
x_{r,i}^{n_{r,i}} \cdot \prod_{k > 0} E_{\beta_k}^{n_k} \,\Big\vert\, n_k,
n_{r,i} \in \N \; \forall \, k, i \,\Big\} \, $  and  $ \, \Big\{\, \prod_{k
\leq 0} F_{\beta_k}^{m_k} \cdot \! \prod_{r \in \N, j \in I_0} y_{s,j}^{m_{s,j}}
\cdot \prod_{k > 0} F_{\beta_k}^{m_k} \,\Big\vert\, m_k, m_{s,j} \in \N \;
\forall \, k, j \,\Big\} \, $  are bases of  $ \uqp $  and  $ \uqm $
respectively which are orthogonal of each other.  More precisely, we have
  $$  \displaylines{
   \pi \left( \prod_{k \leq 0} E_{\beta_k}^{n_k} \cdot \!
\prod_{r \in \N_+, i \in I_0} x_{r,i}^{n_{r,i}} \cdot \prod_{k > 0}
E_{\beta_k}^{n_k} \, , \, \prod_{h \leq 0} F_{\beta_h}^{m_h} \cdot \! \prod_{s
\in \N_+, j \in I_0} y_{s,j}^{m_{s,j}} \cdot \prod_{h > 0} F_{\beta_h}^{m_h}
\right) =  \cr
  = \prod_{\alpha \in \phipre} \delta_{n_\alpha, m_\alpha}
q_\alpha^{\left( n_\alpha \atop 2 \right)} {\, {[n_\alpha]}_{q_\alpha} ! \,
\over \, {\left( q_\alpha^{-1} - q_\alpha \right)}^{n_\alpha} \,} \, \cdot
\prod_{r \in \N, i \in I_0} \delta_{n_{r,i}, m_{r,i}} n_{r,i}! \, \pi {\left(
x_{r,i} \, , y_{r,i} \right)}^{n_{r,i}}  \cr }  $$
\endproclaim

\demo{Proof}  This is just a variation of [Da2], Proposition 10.9, where a
special choice of bases  $ \, \{\, x_{r,i} \mid i \in I_0 \,\} \, $  and  $ \,
\{\, y_{r,j} \mid j \in I_0 \,\} \, $;  but the sole point which is really
necessary is just the orthogonality among such bases: once this is assumed, all
the arguments used in [Da2], \S 10, to build up orthogonal (or even dual) bases
go through as well without change.   $ \square $
\enddemo

\vskip1,7truecm

   \centerline{ \bf  \S \; 3 \,  Integer forms }

\vskip10pt

  {\bf 3.1  The ground ring.} \  As our main goal in  working with quantum
groups is to specialize them at roots of 1, we need integer forms of them
defined over some ring in which  $ (q-\varepsilon) $  be not invertible,
$ \varepsilon $  being any fixed root of 1.  So let  $ {\frak R} $  be the set
of all roots of 1 (in  $ \C $)  whose order is either 1 or an odd number
$ \ell $  with  $ \, g.c.d.(\ell,n+1) = 1 \, $  if  $ \gerg $  is of type
$ A_n \, $,  $ \, \ell \notin 3 \N_+ \, $  if  $ \gerg $  is of type  $ E_6 $
or  $ G_2 \, $;  then set
  $$  \Reps := {\C[q]}_{(q-\varepsilon)} = \left\{\, f \in \C(q) \,\left\vert\,
\hbox{\it  $ f $  has no poles at  $ q = \varepsilon $ } \right. \,\right\} \;
\forall \, \varepsilon \in {\frak R} \, ,  \quad  R := \bigcap_{\varepsilon \in
{\frak R}} \Reps \; ;  $$
in particular,  $ \, \Z\left[q,\qm\right] \subset R \subset \Reps \, $.  Later
on we shall need to inverte the determinants
  $$  \Delta_r := det \left( {\left( {[a_{ij}]}_{q_i^r} \right)}_{\! i,j \in
I_0} \right)  \quad \qquad  \forall \, r \in \N_+ \; ;  $$
of course we have  $ \, \Delta_r \in \Z \left[ q, \qm \right] \, $,  whence in
particular  $ \, \Delta_r \in R \subset \Reps \, $.  Here is the explicit
description of the  $ \Delta_r $'s,  according to the type of  $ \gerg \, $:
  $$  \eqalignno{
   {}\hskip-5cm   A_n :  &  \quad {[n+1]}_{q^r}  &  \qquad  B_n : \quad
{[2]}_{q^{(2n-1)r}}   \hskip3,28cm{}  \cr
   {}\hskip-5cm   C_n :  &  \quad {[2]}_{q^{(n+1)r}}  &  \qquad  D_n : \quad
{[2]}_{q^{(n-1)r}} \cdot {[2]}_{q^r}   \hskip2,48cm{}  \cr
   {}\hskip-5cm   E_6 :  &  \quad {[3]}_{q^r} \cdot \big( {[2]}_{q^{4r}} - 1
\big)  &  \qquad  E_7 : \quad {[2]}_{q^r} \cdot \big( {[2]}_{q^{6r}} - 1
\big)   \hskip2,08cm{}  \cr
   {}\hskip-5cm   E_8 :  &  \quad {[2]}_{q^{8r}} + {[2]}_{q^{6r}} -
{[2]}_{q^{2r}} - 1  &  \qquad  F_4 : \quad {[2]}_{q^{6r}} - 1   \hskip3,34cm{}
\cr
   {}\hskip-5cm   G_2 :  &  \quad {[3]}_{q^r} \cdot \big( {[2]}_{q^{10r}} +
{[2]}_{q^{8r}} - {[2]}_{q^{2r}} - 1 \big) \, .   &  \qquad  {}  \cr }  $$
   \indent   Note that, if  $ \, \varepsilon \in {\frak R} \, $,  then
$ \, {\Delta_r}{\big\vert}_{q=\varepsilon} \neq 0 \, $  for all  $ r \in \N_+
\, $  (by the very definition of  $ \Reps $);  thus
$ {\Delta_r}{\big\vert}_{q=\varepsilon} $  is invertible in each
$ \Reps \, $,  whence also in  $ R \, $.

\vskip7pt

   {\bf 3.2 Lusztig's integer form.} \  Let us define the  $ q $--divided
powers, namely
  $$  F_i^{(\ell)} := {\,F_i^\ell\, \over \, {[\ell]}_{q_i}!\,} \, ,  \qquad
\left[ K_i ; \, c \atop t \right] := \prod_{s=1}^t { \, q_i^{c-s+1} K_i -
q_i^{-(c-s+1)} K_i^{-1} \, \over \, q_i^s - q_i^{-s} \,} \, ,  \qquad
E_i^{(m)} := {\,E_i^m\, \over \,{[m]}_{q_i}!\,}  $$
for all  $ i \in I, \ell, m, t \in \N, c \in \Z $.  Then let  $ \gerUqg $  be
the  $ R $--subalgebra  of  $ \uqg $  generated by
  $$  \left\{\, F_i^{(\ell)}, E_i^{(m)}, K_i^{\pm 1} \,\Big\vert\, i \in I, \,
\ell, m \in \N \,\right\}  \; ;  $$
this is a Hopf  $ R $--subalgebra  of  $ \uqg $  (cf.~[Lu2]).  We set
  $$  \gerUqm := \gerUqg \cap \uqm \, ,  \quad  \gerUqz := \gerUqg \cap
\uqz \, ,  \quad  \gerUqp := \gerUqg \cap \uqp \, ;  $$
then (cf.~[Lu2])  $ \gerUqm $,  resp.~$ \gerUqz $,
resp.~$ \gerUqp $,  is the  $ R $--subalgebra  of  $ \uqg $  generated by
$ \left\{\, F_i^{(\ell)} \,\Big\vert\, \right. $\allowbreak
$ i \in I, \ell \in \N \,\Big\} $,  resp.~$ \left\{\, \left[ K_i ; \, c \atop t
\right] , K_i^{\pm 1} \,\Big\vert\, i \in I_0 \, , t \in \N, c \in \Z
\,\right\} $,  resp.~$ \left\{\, E_i^{(m)} \,\Big\vert\, i \in I, m \in \N
\,\right\} $.
                                                               \par   
   Set
  $$  \gerUqbm := \gerUqg \cap \uqbm \, ,  \quad  \gerUqbp := \gerUqg \cap
\uqbp \, ;  $$
then  $ \gerUqbm $,  resp.~$ \gerUqbp $,  is the  $ R $--subalgebra  of
$ \uqg $  generated by  $ \gerUqm $  and  $ \gerUqz $,  resp.~by  $ \gerUqz $
and  $ \gerUqp $;  moreover, both  $ \gerUqbm $  and  $ \gerUqbp $  are Hopf
subalgebras (over  $ R \, $)  of  $ \gerUqg $.
                                                    \par
  The triangular decompositions in  $ \uqg $  induce similar decompositions
in  $ \gerUqg $,  namely
  $$  \displaylines{
   \gerUqg \cong \gerUqp \otimes \gerUqz \otimes \gerUqm \cong \gerUqm
\otimes \gerUqz \otimes \gerUqp  \cr
   \gerUqbp \cong \gerUqp \otimes \gerUqz \cong \gerUqz \otimes \gerUqp \, ,
\quad  \gerUqbm \cong \gerUqz \otimes \gerUqm \cong \gerUqm \otimes \gerUqz
\, ;  \cr }  $$
which still are given by multiplication.  Finally, the  $ {\Cal B} $--action
on  $ \uqg $  restricts on  $ \gerUqg $.

\vskip4pt

   {\it Remark:} \, Some comments about the "toral part"  $ \gerUqz $  of
$ \gerUqg $  are in order again.  As in the finite case, several choices are
possible: here we made the "minimal one", for we decide simply to take as
$ \gerUqg $  the smallest  $ R $--subalgebra  containing the  $ q $--divided
powers of Chevalley generators  $ E_i $,  $ F_i $  (for all  $ i \, $)  and
stable for the  $ {\Cal B} $--action.

\vskip7pt

   {\bf 3.3 Beck-Kac's integer form.} \  Define renormalized root vectors
by
  $$  \echeck_\alpha := \left( q_\alpha - q_\alpha^{-1} \right) E_\alpha \, ,
\; \quad  \fcheck_\alpha := \left( q_\alpha^{-1} - q_\alpha \right) F_\alpha
\, ,  \; \qquad  \forall \, \alpha \in \phitildep \, ;  $$
(note that  $ \, \fcheck_\alpha = \Omega \left( \echeck_\alpha \right) \, $),
and let  $ \, \calUqg \, $  be the  $ R $--subalgebra of  $ \uqg $  generated
by
  $$  \left\{ \, \fcheck_\alpha, K_\mu, \echeck_\alpha \,\big\vert\, \alpha
\in \Phi_+ \, ,  \mu \in Q_\infty \,\right\} \; .  $$
Then (cf.~[B-K])  $ \, \calUqg \, $  is a Hopf subalgebra of  $ \uqg $,  with
PBW basis over  $ R $
  $$  \Bigg\{\, \prod_{\alpha \in \Phi_+} \echeck_\alpha^{n_\alpha} \cdot
\prod_{i \in I_\infty} K_i^{l_i} \cdot \prod_{\alpha \in \Phi_+}
\fcheck_\alpha^{m_\alpha} \, \Bigg\vert \, l_i \in \Z, n_\alpha, m_\alpha \in
\N, \forall\, i, \alpha \, ; \; \text{almost all} \; \, n_\alpha = m_\alpha = 0
\,\Bigg\}  $$
(here again the monomials are ordered!); the latter is also a  $ \Cq $--basis
of  $ \uqg $,  hence  $ \calUqg $  is an  $ R $--form  of  $ \uqg $;
furthermore, it is a Hopf subalgebra  of  $ \uqg $.  If we define
  $$  \calUqm := \calUqg \cap \uqm \, ,  \quad  \calUqz := \calUqg \cap
\uqz \, ,  \quad  \calUqp := \calUqg \cap \uqp \, ;  $$
then  $ \calUqm $,  resp.~$ \calUqz $,  resp.~$ \calUqp $,  is the
$ R $--subalgebra  of  $ \uqg $  generated by  $ \left\{\, \fcheck_\alpha
\,\big\vert\, \alpha \in \Phi_+ \,\right\} $,  resp.~$ \left\{\, K_i^{\pm 1}
\,\big\vert\, i \in I_\infty \,\right\} $,  resp.~$ \left\{\, \echeck_\alpha
\,\big\vert\, \alpha \in \Phi_+ \,\right\} $.  Letting also
  $$  \calUqbm := \calUqg \cap \uqbm \, ,  \quad  \calUqbp := \calUqg \cap
\uqbp \, ,  $$
we have that  $ \calUqbm $,  resp.~$ \calUqbp $,  is the  $ R $--subalgebra
of  $ \uqg $  generated by  $ \calUqm $  and  $ \calUqz $,  resp.~by
$ \calUqz $  and  $ \calUqp $;  moreover, both  $ \calUqbm $  and
$ \calUqbp $  are Hopf subalgebras (over  $ R \, $)  of  $ \calUqg $.
                                                    \par
  The triangular decompositions in  $ \uqg $  yield similar
decompositions in  $ \calUqg $,  namely
  $$  \displaylines{
   \calUqg \cong \calUqp \otimes \calUqz \otimes \calUqm \cong \calUqm
\otimes \calUqz \otimes \calUqp  \cr
   \calUqbp \cong \calUqp \otimes \calUqz \cong \calUqz \otimes \calUqp \, ,
\quad  \calUqbm \cong \calUqz \otimes \calUqm \cong \calUqm \otimes \calUqz
\cr }  $$
where all the isomorphism are provided by multiplication.
                                                  \par
   It is clear now that PBW bases over  $ R $  also exist for the
previous subalgebras of  $ \calUqg $:  for instance, the set of ordered
monomials in the  $ \, \echeck_\alpha's \, $  is a  $ R $--basis
of  $ \calUqp $;  similarly,  $ \, \left\{\, \prod_{\alpha \in \Phi_+}
\fcheck_\alpha^{m_\alpha} \cdot \prod_{i \in I_\infty } K_i^{t_i}
\,\bigg\vert\, t_i \in \Z, m_\alpha \in \N, \, \text{\ almost all \ } m_\alpha
= 0 \,\right\} \, $  is a  $ R $--basis  of  $ \calUqbm $,  etc.

\vskip4pt

{\it Remark:} \, The previous definitions and results are originally stated in
[B-K] with  $ \Reps $  as ground ring; but since they hold for  {\sl all}  $ \,
\varepsilon \in {\frak R} \, $,  one concludes that they hold over  $ R $  too.

\vskip1,7truecm

   \centerline{ \bf  \S \; 4 \,  The main PBW theorem }

\vskip10pt

   {\bf 4.1  $ q $--divided powers.} \  We extend the notion of  $ q $--divided
powers to the root vectors attached to non-simple roots.  To begin with, set
  $$  \hbox{ $  \eqalign{
    \ehat_\alpha:= E_\alpha \quad \forall \, \alpha \in \phipre \, ,  &
{} \qquad  \ehat_{(r \delta, i)}:= {\, r \, \over \, {[r]}_{q_i} \,} \, E_{(r
\delta,i)} \quad \forall \, (r \delta, i) \in \phitildepim  \cr
    \fhat_\alpha:= F_\alpha \quad \forall \, \alpha \in \phipre \, ,  &
{} \qquad  \fhat_{(r \delta, i)}:= {\, r \, \over \, {[r]}_{q_i} \,} \, F_{(r
\delta,i)} \quad \forall \, (r \delta, i) \in \phitildepim  \cr } $ }
\eqno  (4.1)  $$
note that  $ \, \fhat_\alpha = \Omega \left( \ehat_\alpha \right) $,
$ \ehat_\alpha = \Omega \left( \fhat_\alpha \right) $.  Then define
$ q $--divided  powers of root vectors by
  $$  \hbox{ $  \eqalign{
    E_\alpha^{(k)} \equiv {\ehat_\alpha}^{(k)} := {\, {\ehat_\alpha}^k \, \over
\, {[k]}_{q_\alpha}! \,} \, ,  &  \quad  {\ehat_{(r \delta, i)}}^{(k)} := {\,
{\ehat_{(r \delta, i)}}^k \, \over \, k! \,}  \cr
    F_\alpha^{(k)} \equiv {\fhat_\alpha}^{(k)} := {\, {\fhat_\alpha}^k \, \over
\, {[k]}_{q_\alpha}! \,} \, ,  &  \quad  {\fhat_{(r \delta, i)}}^{(k)} := {\,
{\fhat_{(r \delta, i)}}^k \, \over \, k! \,}  \cr } $ }
        \qquad  \forall \, \alpha \in \phipre \, ,  \; (r \delta, i) \in
\phitildepim \, ,  \; k \in \N   \eqno  (4.2)  $$

\vskip7pt

\proclaim{Lemma 4.2} For all  $ \, \alpha \in \phitildep \, $,  $ \, k \in \N
\, $,  we have  $ \, {\ehat_\alpha}^{(k)} \in \gerUqp \, $,  $ \,
{\fhat_\alpha}^{(k)} \in \gerUqm \, $.
\endproclaim

\demo{Proof}  It is enough to prove the claim for the  $ E $'s.
                                                       \par
  If  $ \, \alpha \in \phipre \, $,  then  $ \, \alpha = w(\alpha_i)
\, $  for some  $ w \in W $  and some  $ i \in I $:  then  $ \, E_\alpha =
T_w \left( E_i \right) \, $  by definition, so that  $ \,
{\ehat_\alpha}^{(k)} = T_w \left( E_i^{(k)} \right) \, $  for all
$ k \in \N $.  But  $ \gerUqg $  is  $ \Cal B $--invariant (cf.~[Lu2]), in
particular  $ \, T_w \left( \gerUqg \right) = \gerUqg \, $,  hence
$ \, {\ehat_\alpha}^{(k)} \in \gerUqg \, $,  q.e.d.
                                                       \par
  If  $ \, \alpha \in \phitildepim \, $,  say  $ \, \alpha = (r \delta, i)
\, $,  we proceed as follows.  Consider the (imaginary) root vectors  $ \,
\edot_{(r \delta, i)} \, $  ($ (r \delta,i) \in \phitildepim $)  defined in \S
2.2: of course they lie in  $ \uqp  \, $;  but even more, it is proved in
[C-P],
\S 5, that in fact it is
  $$  \edot_{(r \delta, i)} \in \gerUqp  \qquad \forall \, (r \delta,i) \in
\phitildepim \, .   \eqno (4.3)  $$
   \indent  Now, the relationship between the  $ \edot_{(h \delta, i)} $'s  and
the  $ E_{(k \delta, i)} $'s  is given by (2.3); reversing that formula we get
  $$  {\, r \, \over \, {[r]}_{q_i} \,} \, E_{(r \delta, i)} = - r q_i^{-r}
\edot_{(r \delta, i)} - \sum_{h=1}^{r-1} q_i^{h-r} {\, h \, \over \,
{[h]}_{q_i} \,} \, E_{(h \delta, i)} \edot_{((r-h) \delta, i)}  \qquad \forall
\, (r \delta,i) \in \phitildepim  $$
which in terms of  $ \ehat_{(t \delta, i)} $'s  reads
  $$  \ehat_{(r \delta, i)} = - r q_i^{-r} \edot_{(r \delta, i)} -
\sum_{h=1}^{r-1} q_i^{h-r} \ehat_{(h \delta, i)} \edot_{((r-h) \delta, i)}
\qquad \forall \, (r \delta,i) \in \phitildepim \, ;  $$
this together with (4.3) gives us by induction
  $$  \ehat_{(r \delta, i)} \in \gerUqp  \qquad \forall \, r \in \N  $$
whence clearly  $ \, {\ehat_{(r \delta, i)}}^{(k)} := \displaystyle{\, 1 \,
\over \, k! \,} \cdot {\ehat_{(r \delta, i)}}^k \in \gerUqp \, $,  q.e.d.
$ \square $
\enddemo

\vskip7pt

\proclaim{Definition 4.3}
  $$  \displaylines{
   \hat{\gerB}_+ := \Bigg\{\, \prod_{\alpha \in \Phi_+}
{\ehat_\alpha}^{(n_\alpha)} \, \Bigg\vert \, n_\alpha \in \N \, \; \forall\, i,
\alpha; \hbox{\rm \ almost all \ } \, n_\alpha = 0 \,\Bigg\} ,  \cr
   \hat{\gerB}_0 := \Bigg\{\, \prod_{i \in I} K_i^{-Ent(t_i/2)} {\left[ K_i ;
\, 0 \atop t_i \right]}_{q_i} \,\Bigg\vert\, t_i \in \N \, \; \forall i \in
I \,\Bigg\} ,  \cr
   \hat{\gerB}_- := \Bigg\{\, \prod_{\alpha \in \Phi_+}
\fhat_\alpha^{(m_\alpha)} \, \Bigg\vert \, m_\alpha \in \N \, \; \forall\, i,
\alpha; \hbox{\rm \ almost all \ } \, m_\alpha = 0 \,\Bigg\} .  \cr }  $$
(where  $ \, Ent(x) := $  {\it integer part of}  $ x \, $).   Then we define
$ \, {\widehat U}_q^+ \, $,  resp.~$ \, {\widehat U}_q^0 \, $,
resp.~$ \, {\widehat U}_q^- \, $,  to be the  $ R $--submodule  of  $ \uqp $,
resp.~$ \uqz $,  resp.~$ \uqm $,  spanned by  $ \hat{\gerB}_+ $,
resp.~$ \hat{\gerB}_0 $,  resp.~$ \hat{\gerB}_- $.
\endproclaim

\vskip4pt

   {\it Remark:} \, By \S 2.3 it is clear that  $ \hat{\gerB}_\star $
($ \star \in \{+,0,-\} $)  is a  $ \Cq $--basis  of  $ U_q^\star $,  hence
$ {\widehat U}_q^\star $  is a  {\sl free}  $ R $--module,  with
$ \hat{\gerB}_\star $  as an  $ R $--basis.

\vskip7pt

\proclaim{Lemma 4.4}  $ \, {\widehat U}_q^0 = \gerUqz \, $;  in other words,
$ \hat{\gerB}_0 $  is a basis of  $ \gerUqz $  over  $ R $.
\endproclaim

\demo{Proof}  The proof is the same as in the finite case: see [Lu1], Theorem
6.7, and references therein.   $ \square $
\enddemo

\vskip7pt

\proclaim{Definition 4.5}  Let  $ K $  be any field, let  $ A $  and  $ B $
be two  $ K $--algebras,  and let  $ \, \pi : A \otimes B \longrightarrow K
\, $  be a perfect  $ K $--bilinear pairing; assume  $ k $  is a subring of
$ K $  and  $ \widehat A $,  resp.~$ \widehat B $,  is an integer form of
$ A $,  resp.~$ B $,  over  $ k $.  We say that  $ \widehat A $  is the
$ k $--dual  of  $ \widehat B $  (with respect to  $ \pi $)  if
   $$  \widehat A = \Big\{\, a \in A \,\Big\vert\, \pi \left( a, \widehat b
\right) \in k  \; \; \forall \, \widehat b \in \widehat B \,\Big\} .  $$
\endproclaim

\vskip7pt

\proclaim{Proposition 4.6}  $ {\widehat U}_q^+ $  and  $ \calUqm $,
resp.~$ {\widehat U}_q^- $  and  $ \calUqp $,  are  $ R $--dual  of each
other with respect to the DRT pairing.
\endproclaim

\demo{Proof}  Consider a positive imaginary root  $ r \delta $,
$ r \in \N_+ \, $;  from \S 2.3 and definitions we have
  $$  \pi \left( \ehat_{(r \delta, i)}, \fcheck_{(r \delta, j)} \right) =
{\big( o(i) o(j) \big)}^r {\, {[r a_{ij}]}_{q_i} \, \over \, {[r]}_{q_i} \,} =
{\big( o(i) o(j) \big)}^r \, {[a_{ij}]}_{q_i^r}  \quad \forall \, i, j \in I_0
\, .  $$
   \indent  Consider the matrix
  $$  M_r := {\left( {\big( o(i) o(j) \big)}^r \, {[a_{ij}]}_{q_i^r}
\right)}_{i,j \in I_0}  $$
we have  $ \, det \big( M_r \big) = \pm \Delta_r \, $,  thus  $ det \big( M_r
\big) $  is invertible in  $ R $,  so the inverse matrix  $ \, M_r^{-1} =
{\big( \mu_{ij} \big)}_{i,j \in I_0} \, $  has all its entries in  $ R \, $;
then define a new basis  $ \big\{\, \,\check{\!\calF}_{(r \delta, j)}
\,\big\vert\, j \in I_0 \,\big\} $  of  $ W_r $  by
  $$  {\,\check{\!\calF}}_{(r \delta,i)} := \sum_{j \in I_0} \mu_{ji}
\fcheck_{(r \delta,j)}  \qquad \quad \forall \, i \in I_0 \, . $$
By construction, we have now
  $$  \pi \left( \ehat_{(r \delta, i)}, \,\check{\!\calF}_{(r \delta, j)}
\right) = \delta_{ij}  \qquad \forall \, i, j \in I_0 \, .  $$
Moreover, since  $ \, \mu_{ij} \in R \, $  (for all  $ i $,  $ j $)  and
imaginary root vectors commute with each other, the set  $ {\Cal B}_- $  of
ordered monomials
  $$  \prod_{h \leq 0} \fcheck_{\beta_h}^{m_h} \cdot \! \prod_{s \in \N, j \in
I_0} \,\check{\!\calF}_{(s \delta, j)}^{m_{(s \delta, j)}} \cdot \prod_{h > 0}
\fcheck_{\beta_h}^{m_h}  $$
is again an  $ R $--basis  of  $ \calUqm $.  The end of the story is that we
can apply Lemma 2.4, hence
  $$  \displaylines {
   \pi \left( \prod_{k \leq 0} \ehat_{\beta_k}^{(n_k)} \cdot \!
\prod_{r \in \N_+, i \in I_0} \ehat_{(r \delta, i)}^{\, (n_{(r
\delta, i)})} \cdot \prod_{k > 0} \ehat_{\beta_k}^{(n_k)} \, , \, \prod_{h
\leq 0} \fcheck_{\beta_h}^{m_h} \cdot \! \prod_{s \in \N_+, j \in I_0}
\,\check{\!\calF}_{(s \delta, j)}^{m_{(s \delta, j)}} \cdot \prod_{h > 0}
\fcheck_{\beta_h}^{m_h} \right) =  \cr
  = q_\alpha^{\sum_{\alpha \in
\phipre} d_\alpha \left( n_\alpha \atop 2 \right)} \prod_{\gamma \in
\phitildep} \delta_{n_\gamma, m_\gamma} \; ;  \cr }  $$
in particular, the last term is  {\sl invertible}  in the ground ring  $ R $.
                                               \par
  Now let  $ \, x \in \uqp \, $:  since  $ \hat{\gerB}_+ $  is a basis of
$ \uqp $  over  $ \Cq $,  we have  $ \, x = \sum_{b \in \hat{\gerB}_+} x_b
b \, $  for some  $ x_b \in \Cq $  (almost all zero).
                                               \par
  Suppose  $ \, \pi(x,y) \in R \, $  for all  $ y \in \calUqm \, $.  Let
$ \, x_{b'} \neq 0 \, $:  by the previous analysis, there exists a (unique)
monomial  $ \beta' $  in the PBW basis  $ {\Cal B}_- $  of  $ \calUqm $  such
that
  $$  \pi \big( x, \beta' \big) = \sum_{b \in {\Cal B}_+} x_b \cdot \pi \big(
b, \beta' \big) = x_{b'} q^s  $$
for some  $ s \in \N $;  since  $ \, \pi \big( x, \beta' \big) \in R \, $  by
hypothesis, we have  $ \, x_{b'} \in R \, $.  We conclude that  $ \, x \in
\widehat{U}_q^+ \, $.  Conversely, if  $ \, x_b \in R \, $  for all
$ \, b \in \hat{\gerB}_+ \, $,  then for any monomial  $ \beta $  in  $ {\Cal
B}_- $  we have
  $$  \pi \left( x, \beta \right) = \sum_{b \in {\Cal B}_+} x_b \cdot \pi
\left( b, \beta \right) = x_{b'} q^s \in R  $$
for a unique  $ \, b' \in {\Cal B}_+ \, $  and some  $ s \in \N $,  whence
$ \, \pi(x,y) \in R \, $  for all  $ y \in \calUqm \, $.
                                               \par
  As  $ \hat{\gerB}_+ $  is an  $ R $--basis  of  $ \widehat{U}_q^+ $,  the
proof is completed.   $ \square $
\enddemo

\vskip7pt

\proclaim{Theorem 4.7}  $ \, \gerU_q^\pm = {\widehat U}_q^\pm \, $;  in other
words,  $ \hat{\gerB}_\pm $  is a basis of  $ \gerU_q^\pm $  over  $ R $.
\endproclaim

\demo{Proof}  By Lemma 4.2 and the definition of  $ \hat{\gerB}_+ $  we get
$ \, \hat{\gerB}_+ \subseteq \gerUqp \, $,  whence  $ \, {\widehat U}_q^+
\subseteq \gerUqp \, $.  On the other hand, thanks to Proposition 4.6 the
converse inclusion will be proved if we show that  $ \, \pi \left( \gerUqp \, ,
\calUqm \right) \subseteq R \, $.  The  "$ - $"  case will follow applying
$ \Omega $.
                                                     \par
  Consider  $ E_i^{(k)} $  ($ i \in I, k \in \N $);  when pairing
$ E_i^{(k)} $  with any monomial  $ {\Cal F} $  in  $ {\Cal B}_- $,  the
formula in Lemma 2.4   --- cf.~also the proof of Proposition 4.6 above ---
gives
  $$  \pi \left( E_i^{(k)} , {\Cal F} \right) = 0  \quad \forall \, {\Cal F}
\neq \fcheck_{\alpha_i}^{(k)} \, ,  \qquad  \pi \left( E_i^{(k)} ,
\fcheck_{\alpha_i}^k \right) =  q^s  $$
for some  $ s \in \N $;  in particular,  $ \, \pi \left( E_i^{(k)} , {\Cal F}
\right) \in R \, $;  therefore  $ \, E_i^{(k)} \in \widehat{U}_q^+ \, $  for
all  $ \, i \in I, k \in \N \, $.  But now remark that  $ \, {\widehat U}_q^+
\, $  is closed for multiplication, so that  $ \, \gerUqp \subseteq
\widehat{U}_q^+ \, $,  q.e.d.  In fact, let  $ \, f, g \in {\widehat U}_q^+
\, $;  since  $ \pi $  is a Hopf pairing, we have
  $$  \pi (f \cdot g, y) = \pi \big( f \otimes g, \Delta(y) \big)  $$
for all  $ \, y \in \calUqm \subseteq \calUqbm \, $;  but  $ \calUqbm $  is a
Hopf  $ R $--algebra,  hence  $ \, \Delta(y) = \sum_{(y)} y_{(1)} \otimes
y_{(2)} \, $  with  $ \, y_{(1)}, y_{(2)} \in \calUqbm \, $,  and so
  $$  \pi (f \cdot g, y) = \pi \big( f \otimes g, \Delta(y) \big) =
\sum_{(y)} \pi \big( f, y_{(1)} \big) \cdot \pi \big( g, y_{(2)} \big)
\in R \, ,   \text{\ q.e.d.}  \; \square  $$
\enddemo

\vskip4pt

   {\it Remark:} \, Notice that, if  $ \, \varepsilon \in {\frak R} \, $  has
order  $ \ell $  greater than 1, then  $ \, {[r]}_{q_i}
{\Big\vert}_{q=\varepsilon} = 0 \, $  for all  $ \, r \in \ell \, \N_+ \, $;
similarly,  $ \, {[n_\alpha]}_{q_\alpha} \! ! {\Big\vert}_{q=\varepsilon} = 0
\, $  for all  $ \, n_\alpha \geq \ell \, $:  this means that both  $ \,
{[r]}_{q_i} \, $  and  $ \, {[m_\alpha]}_{q_\alpha} \! ! \, $  are  {\sl not}
invertible in  $ \Reps \, $,  hence   --- a fortiori ---   no more in
$ R \, $.  Therefore the occurrence of coefficients  $ \, {\left({\, 1 \,
\over \, {[r]}_{q_i}}\right)}^{n_{(r \delta,i)}} \, $  ("hidden inside"  $ \,
\ehat_{(r \delta,i)}^{(n_{(r \delta,i)})} \, $)  and  $ \, {\, 1 \, \over \,
{[n_\alpha]}_{q_\alpha} \! ! \,} \, $  ("hidden inside"  $ \,
\ehat_\alpha^{(n_\alpha)} \, $),  multiplying a monomial of Beck's basis for
$ \, \uqp \, $,  is really significant.

\vskip7pt

   We conclude with our main result (a second PBW theorem is proved in \S 5):

\vskip7pt

\proclaim{First PBW Theorem 4.8}  The sets of ordered monomials  $ \,
\hat{\gerB}_+ \cdot \hat{\gerB}_0 \cdot \hat{\gerB}_- \, $  and  $ \,
\hat{\gerB}_- \cdot \hat{\gerB}_0 \cdot \hat{\gerB}_+ \, $  are bases of  $ \,
\gerUqg \, $  over  $ R \, $.
\endproclaim

\demo{Proof}  Trivial from Lemma 4.4, Theorem 4.7, and triangular decomposition
(cf.~\S 3.2).   $ \square $
\enddemo

\vskip1,7truecm

   \centerline{ \bf  \S \; 5 \,  Generating series, the classical framework,
and beyond }

\vskip10pt

  {\bf 5.1  Changing imaginary root vectors.} \  Let us consider sets of
countably many variables  $ \, {\Bbb X} = {\big\{ X_r \big\}}_{r \in \N} \, $,
with  $ \, X_0 = 1 \, $.  Let such a set  $ {\Bbb X} $  be given: we define a
new set  $ \, {\Bbb Y} = {\big\{ Y_s \big\}}_{s \in \N} \, $  via the following
relation among generating series:
  $$  \sum_{s=0}^\infty Y_s \cdot \zeta^s = \exp \left( \sum_{r=1}^\infty X_r
\cdot \zeta^r \right) \; ;   \eqno (5.1)  $$
(here  $ \zeta $  is any auxiliary symbol) notice this agree with  $ \, Y_0 = 1
\, $.  This is a variation on a classical theme: the  $ Y_s $'s  above are
certain normalizations of the well-known  {\it complete Bell polynomials}
(cf.~[Co], \S 3.3).  We shall shortly denote any "change of variables" given
by (5.1) with  $ \, \Psi \colon {\Bbb X} \strangearrow {\Bbb Y} \, $  or  $ \,
{\Bbb X} {\buildrel \Psi \over \strangearrow} {\Bbb Y} \, $,  and then write
$ \, {\Bbb Y} = \Psi \left( {\Bbb X} \right) \, $.  The converse change, to be
denoted  $ \, \Phi \colon {\Bbb Y} \strangearrow {\Bbb X} \, $  or  $ \, {\Bbb
Y} {\buildrel \Phi \over \strangearrow} {\Bbb X} \, $  or  $ \, {\Bbb X} = \Phi
\left( {\Bbb Y} \right) \, $,  with  $ \, \Phi = \Psi^{-1} \, $,  is given by
  $$  \sum_{r=1}^\infty X_r \cdot \zeta^r = \log \left( \sum_{s=0}^\infty Y_s
\cdot \zeta^s \right) \, .   \eqno (5.2)  $$
Differentiating (with respect to  $ \zeta $)  and comparing the
coefficients of  $ \zeta^r $  on both sides gives
  $$  r Y_r = r X_r + \sum_{s=1}^{r-1} s X_s Y_{r-s} = \sum_{s=1}^r s X_s
Y_{r-s}  \qquad \forall \, r \in \N_+   \eqno (5.3)  $$
which is a recursive rule to compute the  $ Y_r $'s  starting from the
$ X_s $'s,  i.e.~to compute  $ \, {\Bbb Y} = \Psi \left( {\Bbb X} \right) \, $.
Reversing (5.3) gives the rule to compute  $ \, {\Bbb X} = \Phi \left( {\Bbb Y}
\right) \, $,  namely
  $$  r X_r = r Y_r - \sum_{s=1}^{r-1} (r-s) Y_s X_{r-s}  \qquad \forall \, r
\in \N_+ \; .   \eqno (5.4)  $$
Conversely, if  $ {\Bbb X} $  and  $ {\Bbb Y} $  verify (5.3) or (5.4), then
$ \, {\Bbb Y} = \Psi \left( {\Bbb X} \right) \, $  and  $ \, {\Bbb X} = \Phi
\left( {\Bbb Y} \right) \, $.
                                                   \par
   All this applies to our various sets of imaginary root vectors: for every
$ \, i \in I_0 \, $  let
  $$  \displaylines{
   \Etilde_i := {\left\{ \etilde_{(r \delta,i)} \right\}}_{r \in \N} \, ,
\quad  \E_i := {\left\{ E_{(r \delta,i)} \right\}}_{r \in \N} \, ,  \quad
\Edot_i := {\left\{ \edot_{(r \delta,i)} \right\}}_{r \in \N} \, ,  \quad
\Ehat_i := {\left\{ \ehat_{(r \delta,i)} \right\}}_{r \in \N}  \cr
   - \left( q_i - q_i^{-1} \right) \Etilde_i := {\left\{ - \left( q_i -
q_i^{-1} \right) \etilde_{(r \delta,i)} \right\}}_{r \in \N} \, ,  \quad \left(
q_i - q_i^{-1} \right) \E_i := {\left\{ \left( q_i - q_i^{-1} \right) E_{(r
\delta,i)} \right\}}_{r \in \N}  \cr
   - {\, q_i^* \, \over \, \ast \,} \, \Ehat_i := {\left\{ - {\, q_i^r \, \over
\, r \,} \ehat_{(r \delta,i)} \right\}}_{r \in \N} \, ;  \cr }  $$
the relations in \S 2.2 tell us that:  $ \, - \left( q_i - q_i^{-1} \right)
\Etilde_i {\buildrel \Phi \over \strangearrow} \left( q_i - q_i^{-1} \right)
\E_i \, $,  by (2.1), and similarly, by (2.3),  $ \, - {\, q_i^* \, \over \,
\ast \,} \, \Ehat_i {\buildrel \Psi \over \strangearrow} \Edot_i \, $.

\vskip7pt

   {\bf 5.2  The classical setting.} \  Recall that  $ \gerUqg $  is a
quantization of  $ U \big( \ghat \big) $,  that is we have an isomorphism (of
co-Poisson Hopf algebras)
  $$  \gerUqg \Big/ (q-1) \, \gerUqg \,\cong\, U \big( \ghat \big) \, .
\eqno (5.5)  $$
   The construction of the previous section are strictly related with some
classical ones, inside  $ U \big( \ghat \big) $.  Namely, let $ \ghat_\Z $  be
the Lie subalgebra of  $ \ghat $  generated over  $ \Z $  by the Chevalley
generators (assumed to be given) of  $ \ghat \, $;  one can locate a Chevalley
basis of  $ \ghat $,  i.e.~a  $ \Z $--basis  of  $ \ghat_\Z \, $.  Now let
$ \, U_\Z \big( \ghat \big) \, $  be the  $ \Z $--subalgebra  of  $ \, U \big(
\ghat \big) \, $  generated by the divided powers of the Chevalley generators:
this is the classical counterpart of our  $ \gerUqg $,  and the problem of
finding a  $ \Z $--basis  of it is the classical analogue of finding an
$ R $--basis  of  $ \gerUqg $.  The construction of a  $ \Z $--basis  of  $ \,
U_\Z \big( \ghat \big) \, $  resembles that for the finite dimensional
semisimple Lie algebras, but for one point, regarding imaginary root vectors.
                                               \par
  Following [Mi], \S 4, pick out of the Chevalley basis of  $ \ghat $  the
subset  $ \big\{ \, \text{e}_\beta \,\big\vert\, \beta \in \phipre \, \big\} $
of root vectors attached to positive real roots and the subset
$ \big\{ \, \tilde{\text{e}}_{(r \delta,i)} \,\big\vert\, (r \delta, i) \in
\phitildepim \, \big\} $  of root vectors attached to positive imaginary roots
with multiplicity.  For all  $ \, i \in I_0 \, $,  $ \, k \in \N_+ \, $,
define a new set  $ \, \Lam_i^{\langle k \rangle} := {\left\{ \Lambda^{\langle
k \rangle}_{[s;i]} \right\}}_{s \in \N} \, $  by the change of variables
  $$  {\left\{ {\, 1 \, \over \, r \,} \, \tilde{\text{e}}_{(r \cdot k
\delta,i)} \right\}}_{r \in \N} {\buildrel \Psi \over \strangearrow} {\left\{
\Lambda^{\langle k \rangle}_{[s;i]} \right\}}_{s \in \N} \, ,  \quad  \hbox{\
in short \ }  \;  \Lam_i^{\langle k \rangle} := \Psi \left( {\, 1 \, \over \,
\ast \,} \, \tilde{\hbox{\bf e}}_{(k \delta,i)} \right) \; .   \eqno (5.6)  $$
   \indent   Now fix an order on the Chevalley basis of  $ \ghat $,  let  $ \,
b^{(s)} := {\, b^s \, \over \, s! \,} \, $  denote usual divided powers,  let
$ \, \ghat_\Z^+ := \ghat_\Z \bigcap \bigoplus\limits_{\alpha \in \Phi_+}
\ghat \, $,  and set  $ \, U_\Z \left( \ghat^+ \right) :=  U_\Z \left( \ghat
\right) \cap U \left( \ghat^+ \right) \, $  (this is the classical counterpart
of our  $ \gerUqp $).  For any  $ \, {\big( n_\alpha \big)}_{\alpha \in
\phitildep} \in \N^{\phitildep} \, $  such that all but a finite number of
$ n_\alpha $'s  are zero define the  {\sl monomial}  $ \, {\Cal M}^+ \Big(
{\big( n_\alpha \big)}_{\alpha \in \phitildep} \Big) \, $  to be the ordered
product of the elements  $ \, \text{e}_\beta^{(n_\beta)} \, $,  $ \,
\Lambda^{\langle r \rangle}_{\left[ n_{(r \delta,i)}; i \right]} \, $  for all
$ \, \beta \in \phipre \, $,  $ \, (r \delta,i) \in \phitildepim \, $.  The
importance of the new sets of imaginary root vectors  $ \, \Lam_i^{\langle k
\rangle} \, $  (for all  $ i \in I_0 \, $)  lies in the following result, due to
Garland:

\vskip7pt

\proclaim{Theorem 5.3 ([Ga], [Mi])}  The set of all monomials  $ \, {\Cal M}^+
\Big( {\big( n_\alpha \big)}_{\alpha \in \phitildep} \Big) \, $  defined above
is a  $ \Z $--basis  of  $ \, U_\Z \left( \ghat^+ \right) \, $.
$ \square $
\endproclaim

\vskip7pt

   {\bf 5.4  Quantum versus classical imaginary root vectors.} \  Comparison
shows that Theorem 4.8  {\sl is not}  the direct  $ q $--analogue  of its
classical counterpart, Theorem 5.3.
                                                \par
   First, Theorem 4.8 involves divided powers of imaginary root vectors,
whereas Theorem 5.3 involves the  $ \Lambda^{\langle r \rangle}_{\left[ n_{(r
\delta,i)}; i \right]} $'s,  which  {\sl are not}  divided powers of imaginary
root vectors; in this respect, the quantum result is somewhat nicer than the
classical one.
                                                \par
   Second, notice that in the isomorphism (5.5) one has
  $$  E_\alpha {\Big\vert}_{q=1} = \text{e}_\alpha  \quad \forall \, \alpha
\in \phipre \, ,  \qquad  \etilde_{(r \delta,i)} {\Big\vert}_{q=1} =
\tilde{\text{e}}_{(r \delta,i)}  \quad \forall \, (r \delta,i) \in
\phitildepim \; ;   \eqno (5.7)  $$
now, the imaginary root vectors involved in Theorem 5.3 are built from
those in the Cartan basis   --- the  $ \tilde{\text{e}}_{(r \cdot k
\delta, i)} $'s ---   (up to normalization) through a change of variables
"of type  $ \Psi $";  on the other hand, in the quantum setting one starts
from the  $ \etilde_{(r \delta,i)} $'s   --- which are the  $ q $--analogue
of the  $ \tilde{\text{e}}_{(r \cdot k \delta,i)} $'s,  thanks to (5.7) ---
but then the imaginary root vectors occurring in Theorem 4.8 are built from
the  $ \etilde_{(r \delta,i)} $'s  (up to normalizations, for  $ \Ehat _i $
is a renormalization of  $ \E_i \, $)  through a change of variables "of type
$ \Phi $";  in other words, starting from the same point the search for new
suitable imaginary root vectors in the quantum case goes in opposite
direction than in the classical setting.
                                                    \par
   Nevertheless, from Theorem 4.8 we can prove a second PBW Theorem, modeled on
the classical one.  To begin with, applying (5.4) to  $ \, - \left( q_i -
q_i^{-1} \right) \Etilde_i {\buildrel \Phi \over \strangearrow} \left( q_i -
q_i^{-1} \right) \E_i \, $  gives
  $$  r \left( q_i - q_i^{-1} \right) E_{(r \delta,i)} = - r \left( q_i -
q_i^{-1} \right) \etilde_{(r \delta,i)} + \sum_{s=1}^{r-1} (r-s) {\left( q_i -
q_i^{-1} \right)}^2 \etilde_{(s \delta,i)} E_{((r-s) \delta,i)}  $$
which, dividing by  $ \left( q_i - q_i^{-1} \right) $,  implies
  $$  E_{(r \delta,i)}{\Big\vert}_{q=1} = - \etilde_{(r \delta,i)}
{\Big\vert}_{q=1} = - \tilde{\text{e}}_{(r \delta,i)}  \qquad \forall \, (r
\delta,i) \in \phitildepim \; ;   \eqno (5.8)  $$
\noindent   furthermore, we have
  $$  \ehat_{(r \delta,i)}{\Big\vert}_{q=1} = + E_{(r
\delta,i)}{\Big\vert}_{q=1} = - \tilde{\text{e}}_{(r \delta,i)}  \qquad \forall
\, (r \delta,i) \in \phitildepim \, .   \eqno (5.9)  $$
   \indent   Now define a new set of imaginary root vectors  $ \, \Edot_i^{[k]}
= {\left\{ \edot_{[s;i]}^{[k]} \right\}}_{s \in \N} \, $  ($ \, k \in \N_+
\, $)  by
  $$  - {\, q_i^{k \cdot \ast} \, \over \, \ast \,} \Ehat_i(k) := {\left\{ -
{\, q_i^{k \cdot r} \, \over \, r \,} \ehat_{(r \cdot k \delta,i)} \right\}}_{r
\in \N} {\buildrel \Psi \over \strangearrow} {\left\{ \edot_{[s;i]}^{[k]}
\right\}}_{s \in \N} =: \Edot_i^{[k]}   \eqno (5.10)  $$
(in particular,  $ \, \Edot_i^{[1]} = \Edot_i \, $,  by definitions);
explicitly, it is (for all  $ \, r, k \in \N_+ \, $,  $ \, i \in I_0 \, $)
  $$  \edot_{[0;i]}^{[k]} := 1 \, ,  \quad  \edot_{[r;i]}^{[k]} :=
- {\, 1 \, \over \, r \,} \sum_{s=1}^r q_i^{k \cdot s} {\, s \cdot k \, \over
\, {[s \cdot k]}_{q_i} \,} E_{(s \cdot k \delta, i)} \edot_{[(r-s);i]}^{[k]}
\quad \forall \, (r \delta,i) \in \phitildepim \; .   \eqno (5.11)  $$
Since  $ \, \Ehat_i \subset \gerUqp \, $,  we have also  $ \, \Ehat_i(k)
\subset \gerUqp \, $,  whence  $ \, \Edot_i^{[k]} \subset \gerUqp \, $;
moreover, from (5.6--10) we have  $ \, \Edot_i^{[k]}{\Big\vert}_{q=1} =
\Lam_i^{\langle k \rangle} \, $:  in other words, the
$ \edot_{[r;i]}^{[k]} $'s  are  $ q $--analogue  of the
$ \Lambda^{\langle k \rangle}_{[r;i]} $'s.
                                                       \par
   Now we are ready for our second PBW theorem, which strictly mimicks the
classical one.  First a further definition: for any  $ \, {\big( n_\alpha
\big)}_{\alpha \in \phitildep} \in \N^{\phitildep} \, $  such that all but a
finite number of  $ n_\alpha $'s  are zero define the  {\sl monomial}  $ \,
{\Cal M}^+_q \Big( {\big( n_\alpha \big)}_{\alpha \in \phitildep} \Big) \, $
to be the ordered product   --- with respect to the order of  $ \phitildep $
fixed in \S 2.1 ---   of the elements  $ \, E_\beta^{(n_\beta)} \, $,
$ \, \edot^{[r]}_{\left[ n_{(r \delta,i)}; i \right]} \, $  for all  $ \, \beta
\in \phipre \, $,  $ \, (r \delta,i) \in \phitildepim \, $;  similarly define
$ \, {\Cal M}^-_q \Big( {\big( n_\alpha \big)}_{\alpha \in \phitildep} \Big)
\, $  by means of negative root vectors (with  $ \, \fdot_{[r;i]}^{[k]} :=
\Omega \left( \edot_{[r;i]}^{[k]} \right) \, $  of course), and call
$ \, \dot{\gerB}_\pm^{[\ast]} \, $  the set of all monomials  $ {\Cal M}^\pm_q
\Big( {\big( n_\alpha \big)}_{\alpha \in \phitildep} \Big) $.

\vskip7pt

\proclaim{Second PBW Theorem 5.5}    \phantom{$\vert$}
                                        \hfill\break
   \indent   (a) \  $ \dot{\gerB}_+^{[\ast]} \, $,  resp.~$ \,
\dot{\gerB}_-^{[\ast]} \, $,  is an  $ R $--basis  of  $ \, \gerUqp \, $,
resp.~$ \, \gerUqm \, $.
                                      \hfill\break
   \indent   (b) \  The sets  $ \, \dot{\gerB}_+^{[\ast]} \cdot \hat{\gerB}_0
\cdot \dot{\gerB}_-^{[\ast]} \, $  and  $ \, \dot{\gerB}_-^{[\ast]} \cdot
\hat{\gerB}_0 \cdot \dot{\gerB}_+^{[\ast]} \, $  are  $ R $--bases  of
$ \gerUqg $.
\endproclaim

\demo{Proof}  Of course  {\it (b)}  follows from  {\it (a)}  and triangular
decomposition.
                                                           \par
   As for  {\it (a)},  we proceed as follows.  Inside  $ U \left( \ghat^+
\right) $  we have the basis provided by Theorem 5.3 and the usual basis of
ordered monomials in the positive root vectors (of the given Chevalley basis):
in particular, the first basis involve imaginary (positive) root vectors
$ \tilde{\text{e}}_{(r \cdot k \delta,i)} $,  whereas the second one involve
the  $ \Lambda^{\langle k \rangle}_{[r;i]} $,  and the relationship among the
two kinds of vectors is given by  $ \, {\, 1 \, \over \, \ast \,} \,
\tilde{\hbox{\bf e}}_{(k \delta,i)} {\buildrel \Psi \over \strangearrow}
\Lam_i^{\langle k \rangle} \, $.  As they are both bases of  $ U \left( \ghat^+
\right) $  over  $ \C $,  the matrix of basis change has entries in  $ \C $
(even more, in  $ \Q $).
                                                           \par
   Now, the situation in the quantum context is exactly the same; consider the
set  $ \hat{\gerB}'_+ $  defined like  $ \hat{\gerB}_+ $  but for the fact that
every  $ \ehat_{(z \delta,i)} $  is replaced by  $ q_i^z \ehat_{(z \delta,
i)} $:  of course  $ \hat{\gerB}'_+ $  is again an  $ R $--basis  of
$ \gerUqp $  (thanks to Theorem 4.7).  Then  $ \hat{\gerB}'_+ $  and
$ \dot{\gerB}_+^{[\ast]} $  differ only in the fact that they are built up in
the same way but for the use of  $ \, {\, 1 \, \over \, \ast \,} \left( {-
q_i^{k \cdot \ast} \Ehat_i(k)} \right) \, $  in the first case and  $ \,
\Edot_i^{[k]} \, $  in the second one, and their link is again given by  $ \,
{\, 1 \, \over \, \ast \,} \left( {- q_i^{k \cdot \ast} \Ehat_i(k)} \right)
{\buildrel \Psi \over \strangearrow} \Edot_i^{[k]} \, $.  Therefore in force of
the previous analysis these two sets must have the same  $ \C $--span  (even
more, the same  $ \Q $--span),  whence the claim follows.   $ \square $
\enddemo

\vskip7pt

   {\bf 5.6  Integrality questions and beyond.} \  Our constructions and
results assume  $ R $  as ground ring; on the other hand, Lusztig's original
definition of  $ \gerUqg $  takes  $ \, \Z \left[ q, \qm \right] \, $  as
ground ring: thus the question arises of extending our results to such a
setting.
                                               \par
   Conjecturally, the definition of Beck-Kac's integer forms and their
properties should hold as well over
                                              \par
   \indent \indent   {\it (a)} \, the smaller ring  $ \, \C \left[ q, \qm
\right] \, $;
                                              \par
   \indent \indent   {\it (b)} \, the even smaller ring  $ \, \Q \left[ q, \qm
\right] \, $;
                                              \par
   \indent \indent   {\it (c)} \, the smallest ring  $ \, \Z \left[ q, \qm
\right] \, $.
                                              \par
   If  $ W $  is anyone of the above rings, let  $ {\Cal S}_W $  be the
multiplicative part of  $ W $  generated by  $ \{\, \Delta_r \mid r \in \N_+
\,\} \, $,  and set  $ \,W_{\Cal S} := {\Cal S}^{-1}_W W \, $.
                                              \par
   In case  {\it (a)}  or  {\it (b)}  would be true, it is immediate to check
that all our arguments would go through as well if the ground ring is  $ \, {\C
\left[ q, \qm \right]}_{\Cal S} \, $  in case  {\it (a)}  or if it is  $ \, {\Q
\left[ q, \qm \right]}_{\Cal S} \, $  in case  {\it (b)}.  Case  {\it (c)}  is
more tricky.  First of all, we remark that an argument used in the proof of
Lemma 4.2 can be refined: namely, formula (19) in [C-P] shows that in fact
$ \edot_{(r \delta,i)} $  belongs to the Lusztig's form defined over  $ \Z
\left[ q, \qm \right] $;  therefore, the same steps in the cited proof tell the
same is true for  $ \ehat_{(r \delta,i)} $  too, but this does no longer imply
that the same holds for its divided powers  $ \ehat_{(r \delta,i)}^{(k)} $
($ k \in \N $).  In fact for sure Theorem 4.8  {\sl cannot hold}  over  $ {\Z
\left[ q, \qm \right]}_{\Cal S} \, $,  because otherwise by specialization at
$ \, q=1 \, $  it would imply   --- since  $ \, E_\beta^{(k)} {\Big\vert}_{q=1}
\! = \text{e}_\beta^{(k)} $,  $ \, \ehat_\gamma^{(k)} {\Big\vert}_{q=1} \! =
\tilde{\text{e}}_\gamma^{(k)} $,  for all  $ \, \beta \in \phipre \, $,
$ \, \gamma \in \phitildepim \, $  ---   that the set of ordered products of
divided powers of root vectors is a  $ \Z $--basis  of  $ U_\Z \left( \ghat^+
\right) $,  which is  {\sl false}.
                                            \par
   Things change for Theorem 5.5: its statement is really a  $ q $--analogue
of the classical one, which relies on Theorem 5.3.  This lead us to make the
following

\vskip7pt

\proclaim{Conjecture A}    \phantom{$\vert$}
                                         \hfill\break
   \indent   (a) \  $ \dot{\gerB}_+^{[\ast]} \, $,  resp.~$ \,
\dot{\gerB}_-^{[\ast]} \, $,  is a  $ \Z \left[ q, \qm \right] $--basis  of the
Lusztig's form of  $ \uqp $,  resp.~$ \uqm $,  defined  over  $ \Z \left[ q,
\qm \right] \, $.
                                      \hfill\break
   \indent   (b) \  The sets  $ \, \dot{\gerB}_+^{[\ast]} \cdot \hat{\gerB}_0
\cdot \dot{\gerB}_-^{[\ast]} \, $  and  $ \, \dot{\gerB}_-^{[\ast]} \cdot
\hat{\gerB}_0 \cdot \dot{\gerB}_+^{[\ast]} \, $  are  $ \Z \left[ q, \qm
\right] $--bases  of the Lusztig's form of  $ \uqg $  defined  over
$ \Z \left[ q, \qm \right] \, $.
\endproclaim

\vskip7pt

   Furthermore, basing on the classical construction we can sketch an
alternative approach.
                                     \par
   Recall that in the classical setting to achieve a PBW theorem over  $ \Z $
one begins by the change of variables (5.6), i.e.  $ \, {\left\{ {\, 1 \, \over
\, r \,} \, \tilde{\text{e}}_{(r \cdot k \delta,i)} \right\}}_{r \in \N}
{\buildrel \Psi \over \strangearrow} {\left\{ \Lambda^{\langle k
\rangle}_{[r;i]} \right\}}_{r \in \N} \, $:  such a change is described either
by a relation between generating series   --- namely any one of (5.1), (5.2)
---   or by a recursive formula   --- any one of (5.3), (5.4).  In fact, in the
original approach, to be found in [Ga], the  $ \, \Lambda^{\langle k
\rangle}_{[r;i]} $'s  are defined by means of a derivation operator and a
multiplication operator (one for each  $ i \, $)  acting on the  $ \,
\tilde{\text{e}}_{(r \cdot k \delta,i)} $'s;  and then one finds that the  $ \,
\Lambda^{\langle k \rangle}_{[r;i]} $'s  verify the cited recursive formulas.
As a first attempt, one might try to "quantize" this method: exactly this is
partially done in [C-P].  There the  $ q $--analogue  of the  $ \,
\Lambda^{\langle k \rangle}_{[r;i]} $'s,  namely the  $ \, \edot_{(r \delta,i)}
\, $'s  ($ \, \equiv \edot^{[1]}_{[r;i]} \, $),  are defined by means of (2.2),
which is nothing but a straight quantization of a recursive formula of type
(5.3): after this the properties of the  $ \, \edot^{[1]}_{[r;i]} \, $'s  are
obtained by working with a derivation operator and a multiplication operator
which are the  $ q $--analogue  of those used by Garland.  To complete the
task, one should generalize this by defining new vectors  $ \, \edot^{\langle k
\rangle}_{[r;i]} \, $  for all  $ \, k \in \N_+ \, $.  First rewrite the right
hand side formula in (2.2) as
  $$  {(r)}_{q_i^{-2}} \edot_{(r \delta,i)} := \sum_{s=1}^r {(s)}_{q_i^{-2}}
\left( {\, q_i \, \over \, {(s)}_{q_i^{-2}} \,} \, \etilde_{(s \delta, i)}
\right) \edot_{((r-s) \delta,i)}  \qquad \forall \, (r \delta,i) \in
\phitildepim  $$
where we used the well-known notation  $ \, {(z)}_q := {\, q^z - 1 \, \over \,
q - 1 \,} \, $;  this gives (2.2) the shape of a relation of type (5.3).  Now
generalize this defining  $ \, {\left\{ \edot^{\langle k \rangle}_{[r;i]}
\right\}}_{r \in \N} \, $  by means of
  $$  {(r)}_{q_i^{-2}} \edot^{\langle k \rangle}_{[r;i]} := \sum_{s=1}^r
{(s)}_{q_i^{-2}} \left( {\, q_i \, \over \, {(s)}_{q_i^{-2}} \,} \, \etilde_{(s
\cdot k \delta, i)} \right) \edot^{\langle k \rangle}_{[(r-s);i]}  \qquad
\forall \, i \in I_0 \, , \, r, k \in \N_+   \eqno (5.12)  $$
  $$  \hbox{i.e.}  \;\;\;\; \hfill  \edot^{\langle k \rangle}_{[0;i]} := 1,  \;
 \edot^{\langle k \rangle}_{[r;i]} := {\, 1 \, \over \, {(r)}_{q_i^{-2}} \,}
\sum_{s=1}^r {(s)}_{q_i^{-2}} \left( {\, q_i \, \over \, {(s)}_{q_i^{-2}} \,}
\, \etilde_{(s \cdot k \delta, i)} \right) \edot^{\langle k \rangle}_{[(r-s);i]}
 \;\;\; \forall \, i \! \in \! I_0, r, k \! \in \! \N_+  $$
   \indent   Since  $ \, \etilde_{(s \delta,i)} {\Big\vert}_{q=1} =
\tilde{\text{e}}_{(s \delta,i)} \, $  and  $ \, {(r)}_{q_i^{-2}}
{\Big\vert}_{q=1} = r \, $,  definitions give  $ \, \edot^{\langle k
\rangle}_{[r;i]} {\Big\vert}_{q=1} = \Lambda^{\langle k \rangle}_{[r;i]} \, $,
i.e.~the  $ \, \edot^{\langle k \rangle}_{[r;i]} \, $'s  are  $ q $--analogue
of the  $ \, \Lambda^{\langle k \rangle}_{[r;i]} \, $'s;  the same is also true
for the  $ \, \edot^{[k]}_{[r;i]} \, $'s  (cf.~\S 5.4) and in fact  $ \,
\edot^{\langle 1 \rangle}_{[r;i]} = \edot_{(r \delta,i)} = \edot^{[1]}_{[r;i]}
\, $,  but in general it is  $ \, \edot^{\langle k \rangle}_{[r;i]} \neq
\edot^{[k]}_{[r;i]} \, $.  Therefore, if we define the set
$ \dot{\gerB}_+^{\langle \ast \rangle} $  just like the
$ \dot{\gerB}_+^{[\ast]} $  but for using the  $ \, \edot^{\langle k
\rangle}_{[r;i]} \, $'s  instead of the  $ \, \edot^{[k]}_{[r;i]} \, $'s,  and
similarly for  $ \dot{\gerB}_-^{\langle \ast \rangle} $,  then the following
statement should be the most proper  $ q $--analogue  of the classical PBW
theorem over  $ \Z $ (relying on Theorem 5.3):

\vskip7pt

\proclaim{Conjecture B}    \phantom{$\vert$}
                                         \hfill\break
   \indent   (a) \  $ \dot{\gerB}_+^{\langle \ast \rangle} \, $,  resp.~$ \,
\dot{\gerB}_-^{\langle \ast \rangle} \, $,  is a  $ \Z \left[ q, \qm
\right] $--basis  of the Lusztig's form of  $ \uqp $,  resp.~$ \uqm $,  defined
over  $ \Z \left[ q, \qm \right] \, $.
                                      \hfill\break
   \indent   (b) \  The sets  $ \, \dot{\gerB}_+^{\langle \ast \rangle} \cdot
\hat{\gerB}_0 \cdot \dot{\gerB}_-^{\langle \ast \rangle} \, $  and  $ \,
\dot{\gerB}_-^{\langle \ast \rangle} \cdot \hat{\gerB}_0 \cdot
\dot{\gerB}_+^{\langle \ast \rangle} \, $  are  $ \Z \left[ q, \qm
\right] $--bases  of the Lusztig's form of  $ \uqg $  defined  over
$ \Z \left[ q, \qm \right] \, $.
\endproclaim

\vskip7pt

   A way to prove the previous conjecture might pass through the following
alternative approach: quantizing the method in [Mi] instead of that in [Ga].
In fact in [Mi], one defines the  $ \, \Lambda^{\langle k \rangle}_{[r;i]} $'s
via (5.1)   --- or (5.2), it is equivalent ---   then establishes several
relations involving generating series, exponentials and logarithms, and finally
from these recovers Theorem 5.3 (and the most general result for the whole
$ U_\Z \big( \ghat \big) \, $).  So in order to quantize this track, one should
try to directly quantize (5.1), (5.2), instead of (5.3), (5.4), then find
$ q $--analogues  of the classical relations among generating series and the
other functions involved in the "quantized (5.1), (5.2)", and from these
achieve the claim of Conjecture B.  Here we sketch some ideas which might lead
to complete at least the first step.
                                                            \par
   Starting from (5.2), one gets (5.3) by
  $$  {\Bbb X}(\zeta) = \log \big( {\Bbb Y}(\zeta) \big)  \; \Longrightarrow
\;  {\Cal D} \big( {\Bbb X}(\zeta) \big) = {\, {\Cal D} \big( {\Bbb
Y}(\zeta) \big) \, \over \, {\Bbb Y}(\zeta) \,}  \; \Longrightarrow \;  {\Cal
D} \big( {\Bbb Y}(\zeta) \big) = {\Cal D} \big( {\Bbb X}(\zeta) \big) \cdot
{\Bbb Y}  $$
(where  $ \, {\Bbb X}(\zeta) := \sum_{r=0}^\infty X_r \cdot \zeta^r \, $,  $ \,
{\Bbb Y}(\zeta) := \sum_{s=0}^\infty Y_s \cdot \zeta^s \, $,  and
$ {\Cal D} $  denotes derivation with respect to  $ \zeta $)  and
comparison of coefficients of  $ \zeta^t \, (t \in \N) $  on both sides.  Then
the coefficient  $ r $,  resp.~$ s $,  in the left, resp.~right, hand side of
(5.2) just arises from  $ \, {\Cal D} \left( \zeta^r \right) \, $,  resp.~$ \,
{\Cal D} \left( \zeta^s \right) \, $.  Similarly, in (5.12) there are
coefficients  $ {(r)}_{q_i^{-2}} $  and  $ {(s)}_{q_i^{-2}} $  which
                         should come out in a\break
\noindent   similar way: thus, consider the skew-derivation  $ \, {\Cal D}_q
: A((\zeta)) \rightarrow A((\zeta)) \, $  (here  $ A $  is any algebra
containing the variables involved, namely the  $ X_r $'s,  the  $ Y_s $'s,  and
so on) defined by
  $$   {\Cal D}_q : f(\zeta) \mapsto \big( {\Cal D}_q (f) \big) (\zeta) := {\,
f(q \zeta) - f(\zeta) \, \over \, (q-1) \zeta \,}   \qquad \forall \, f(\zeta)
\in A((\zeta)) \; ;  $$
then we have  $ \, {\Cal D}_q \big( \zeta^n \big) = {(n)}_q \zeta^{n-1} \, (n
\in \N) \, $,  hence  $ \, {\Cal D}_{q_i^{-2}} \, $  would fulfill our
requirements.  Second, a suitable function  $ {\lambda o \gamma}_q $
($ q $--analogue  of the logarithm) should be found such that  $ \, {\Cal
D}_{q_i^{-2}} \left( {\lambda o \gamma}_{q_i^{-2}} \left( \Edot^{\langle k
\rangle}_i (\zeta) \right) \right) = {\, {\Cal D}_{q_i^{-2}} \left(
\Edot^{\langle k \rangle}_i (\zeta) \right) \, \over \, \Edot^{\langle k
\rangle}_i (\zeta) \,} \, $  (where  $ \, \Edot^{\langle k \rangle}_i (\zeta)
:= \sum_{r=0}^\infty \edot^{\langle k \rangle}_{[r;i]} \cdot \zeta^r \, $);  to
this end, it should be useful the relation among generating series
  $$  \Edot^{\langle k \rangle}_i \left( q_i^{-2} \zeta \right) = \left( 1 -
\left( q_i - q_i^{-1} \right) \sum_{r=1}^\infty \etilde_{(r \cdot k,i)} \cdot
\zeta^r \right) \cdot \Edot^{\langle k \rangle}_i (\zeta)  $$
which arises from a reformulation of (5.12).

\vskip1,1truecm

\vskip0,7truecm

{}

\enddocument

\vskip1,1truecm

\Refs
\endRefs

\vskip7pt

\smallrm

[Be1] \  J.~Beck,  {\smallit Braid group action and quantum affine
algebras\/},  Commun.~Math.~Phys.~{\smallbf 165} (1994), 555--568.

\vskip4pt

[Be2] \  J.~Beck,  {\smallit Convex bases of PBW type for quantum affine
algebras\/},  Commun.~Math.~Phys.~{\smallbf 165} (1994), 193--199.

\vskip4pt

[B-K] \  J.~Beck, V.~G.~Kac,  {\smallit Finite dimensional representations of
quantum affine algebras at roots of 1\/},  J.~Amer.~Math.~Soc.~{\smallbf 9}
(1996), 391--423.

\vskip4pt

[Bo] \  N.~Bourbaki,  {\smallit Groupes et alg\`ebres de Lie\/},  Chapitres
4--6,
Hermann, Paris, 1968.

\vskip4pt

[Co] \  L.~Comtet,  {\smallit Advanced Combinatorics\/},  D.~Reidel Publishing
Company, Dordrecht--Holland/Boston--U.S.A., 1974.

\vskip4pt

[C-P] \  V.~Chari, A.~Pressley,  {\smallit Quantum affine algebras at roots
of unity\/},  Preprint q-alg/9609031.

\vskip4pt

[Da1] \  I.~Damiani,  {\smallit The highest coefficient of  $ H_\eta $  and
the center at odd roots of 1 for untwisted affine quantum algebras\/},
J.~Algebra {\smallbf 186} (1996), 736--780.

\vskip4pt

[Da2] \  I.~Damiani,  {\smallit La  R--matrice  pour les alg\`e{}bres
quantiques de type affine non tordu\/},  Preprint.

\vskip4pt

[Dr] \  V.~G.~Drinfeld,  {\smallit Quantum groups\/},  Proc.~ICM Berkeley 1
(1986), 789--820.

\vskip4pt

[Ga] \  H.~Garland,  {\smallit The arithmetic theory of loop algebras\/},
J.~Algebra~{\smallbf 53} (1978), 480--551.

\vskip4pt

[Ka] \  V.~G.~Kac,  {\smallit Infinite Dimensional Lie Algebras\/},
Birkh\"auser, Boston, 1983.

\vskip4pt

[Lu1] \  G.~Lusztig,  {\smallit Quantum groups at roots of 1\/},
Geom.~Dedicata~{\smallbf 35} (1990), 89--113.

\vskip4pt

[Lu2] \  G.~Lusztig,  {\smallit Introduction to quantum groups\/},  Progress in
Mathematics {\smallbf 110}, Birkh\"auser, Boston, 1993.

\vskip4pt

[Mi] \  D.~Mitzman,  {\smallit Integral bases for affine Lie algebras and their
universal enveloping algebras\/},  Cont. Math.~{\smallbf 40} (1985).

\vskip4pt

[Ta] \  T.~Tanisaki,  {\smallit Killing forms, Harish-Chandra Isomorphisms, and
Universal R-Matrices for Quantum Algebras\/},  Internat.~J.~Modern Phys.~A
{\smallbf 7}, Suppl.~1B (1992), 941--961.

\vskip0,7truecm

{}


\Refs
  \widestnumber\key {B-K}

\vskip5pt

\ref
  \key  Be1   \by  J. Beck
  \paper  Braid group action and quantum affine algebras
  \jour  Commun. Math. Phys.   \vol  165   \yr  1994   \pages  555--568
\endref

\ref
  \key  Be2   \by  J. Beck
  \paper  Convex bases of PBW type for quantum affine algebras
  \jour  Commun. Math. Phys.   \vol  165   \yr  1994   \pages  193--199
\endref

\ref
  \key  B-K   \by  J. Beck, V. G. Kac
  \paper  Finite dimensional representations of quantum affine algebras
at roots of 1
  \jour  J. Amer. Math. Soc.   \vol  9   \yr  1996   \pages  391--423
\endref

\ref
  \key  Bo   \by  N. Bourbaki
  \book  Groupes et alg\`ebres de Lie, Chapitres 4--6
  \publ  Hermann   \publaddr  Paris   \yr  1968
\endref

\ref
  \key  Co   \by  L. Comtet
  \book  Advanced Combinatorics
  \publ  D. Reidel Publishing Company   \publaddr
Dordrecht--Holland/Boston--U.S.A.   \yr  1974
\endref

\ref
  \key  C-P   \by  V. Chari, A. Pressley
  \paper  Quantum affine algebras at roots of unity
  \jour  preprint  math.QA/9609031
\endref

\ref
  \key  Da1   \by  I. Damiani
  \paper  The highest coefficient of  $ H_\eta $  and the center at odd
roots of 1 for untwisted affine quantum algebras
  \jour  J. Algebra   \vol  186   \yr  1996   \pages  736--780
\endref

\ref
  \key  Da2   \by  I. Damiani
  \paper  La  R--matrice  pour les alg\`e{}bres quantiques de type affine
non tordu
  \toappear  \ in Ann. Scient. \'Ec. Norm. Sup.,  4$^e$  s\'erie, t. 31, 1998
\endref

\ref
  \key  Dr   \by  V. G. Drinfeld
  \paper  Quantum groups
  \jour  Proc.~ICM Berkeley   \vol  1   \yr  1986   \pages  789--820
\endref

\ref
  \key  Ga   \by  H. Garland
  \paper  The arithmetic theory of loop algebras
  \jour J. Algebra    \vol  53   \yr  1978   \pages  480--551
\endref

\ref
  \key  Ka   \by  V. G. Kac
  \book  Infinite Dimensional Lie Algebras
  \publ  Birkh\"auser   \publaddr  Boston   \yr  1983
\endref

\ref
  \key  Lu1   \by  G. Lusztig
  \paper  Quantum groups at roots of 1
  \jour  Geom. Dedicata   \vol  35   \yr  1990   \pages  89--113
\endref

\ref
  \key  Lu2   \by  G. Lusztig
  \paper  Introduction to quantum groups
  \book  Progress in Mathematics  \vol  110
  \publ  Birkha\"user   \publaddr  Boston   \yr  1993
\endref

\ref
  \key  Mi   \by  D. Mitzman
  \paper  Integral bases for affine Lie algebras and their
universal enveloping algebras
  \jour  Cont. Math.   \vol  40  \yr  1985
\endref

\ref
  \key  Ta   \by  T. Tanisaki
  \paper  Killing forms, Harish-Chandra Isomorphisms, and
Universal  $R$--Matrices  for Quantum Algebras
  \jour  Internat. J. Modern Phys. A  \vol 7, Suppl. 1B   \yr  1992
\pages  941--961
\endref

\endRefs
\enddocument